\RequirePackage{fix-cm}
\documentclass[times]{article}          
\usepackage{graphicx}
\usepackage{amsmath}

\usepackage{stmaryrd}

\usepackage{caption}
\usepackage{subcaption}
\usepackage[colorlinks,citecolor=blue,linkcolor=blue,urlcolor=blue]{hyperref}

\usepackage{doi}
\usepackage[authoryear]{natbib}
\setcitestyle{authoryear,open={(},close={)}}

\usepackage{float}

\title{Hole Seeding in Level Set Topology Optimization via Density Fields}

\author{Jorge L. Barrera$^1$, Markus J. Geiss$^2$, Kurt Maute $^1$ \\[12pt]
$^1$Ann and H. J. Smead Department of Aerospace Engineering Sciences,\\ 
University of Colorado at Boulder, 3775 Discovery Dr, Boulder, CO 80303, USA\\
$^2$OHB System AG, Manfred-Fuchs-Strasse 1, 82234 Wessling, Germany \\ 
Corresponding author: maute@colorado.edu}

\begin{document}




\maketitle

\begin{abstract}
Two approaches that use a density field for seeding holes in level set topology optimization are proposed. In these approaches, the level set field describes the material-void interface while the density field describes the material distribution within the material phase. Both fields are optimized simultaneously by coupling them through either a single abstract design variable field or a penalty term introduced into the objective function. These approaches eliminate drawbacks of level set topology optimization methods that rely on seeding the initial design domain with a large number of holes. Instead, the proposed approaches insert holes during the optimization process where beneficial. The dependency of the optimization results on the initial hole pattern is reduced, and the computational costs are lowered by keeping the number of elements intersected by the material interface at a minimum.
In comparison to level set methods that use topological derivatives to seed small holes at distinct steps in the optimization process, the proposed approaches introduce holes continuously during the optimization process, with the hole size and shape being optimized for the particular design problem. 
The proposed approaches are studied using the extended finite element method for spatial discretization, and the solid isotropic material with penalization for material interpolation using fictitious densities. Their robustness with respect to algorithmic parameters, dependency on the density penalization, and performance are examined through 2D and 3D benchmark linear elastic numerical examples, and a geometrically complex mass minimization with stress constraint design problem.
\end{abstract}

\section{Introduction} \label{sec:intro}

Topology optimization (TO) via level set (LS) methods typically relies on shape sensitivities along material interfaces to evolve the design. Considering material-void problems, holes cannot nucleate, only merge, split or disappear.
To facilitate topological changes in the optimization process, either the initial design domain is seeded with an array of holes or holes are introduced at distinct stages of the optimization process. 
In this paper, two alternative hole seeding approaches that nucleate holes during the design process are proposed.

\begin{figure*}[h]
	\center
	\includegraphics[width=0.95\linewidth]{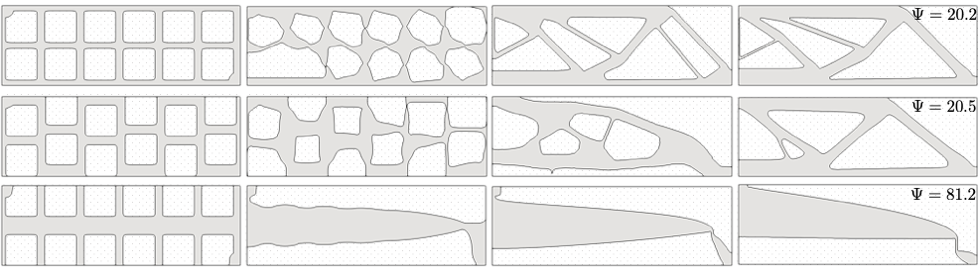}
	\caption{Stages of the optimization process for a 2D Beam LS TO problem that minimizes the strain energy ($\Psi$) using different initial holes seeding patterns.}
	\label{fig:holesArrangementDependency}
\end{figure*}

In absence of any additional mechanism that introduces holes in the course of the optimization process, appropriate initial hole seeding is critical to avoid suboptimal designs (\cite{van2013level}).
In cases where the geometry and the associated mechanical properties of a particular hole pattern violate design constraints (e.g., mass, stress, stiffness), often the majority of holes merge in the early stages of the evolution of the design. 
Premature hole merging can restrict the final design to simplistic geometries and, in certain cases, trigger the appearance of large disconnected subdomains.
The latter can result in ill-posed problems and may compromise the stability of the TO process.
Even small variations in seeding parameters (i.e., number, size, shape or arrangement of holes) can degenerate the evolution of the design and severely affect its performance (\cite{wang2007hole}).
An example of this scenario is shown in Fig. \ref{fig:holesArrangementDependency} for a strain energy minimization problem subject to a mass constraint (\cite{bendsoe2004topology}). 
Stages of the design process with different initial seeding patterns are shown to demonstrate that early merging of holes (bottom row) can result in final designs with poor performance ($>400\%$ difference).

\begin{figure}[h]
	\center
	\includegraphics[width=0.5\linewidth]{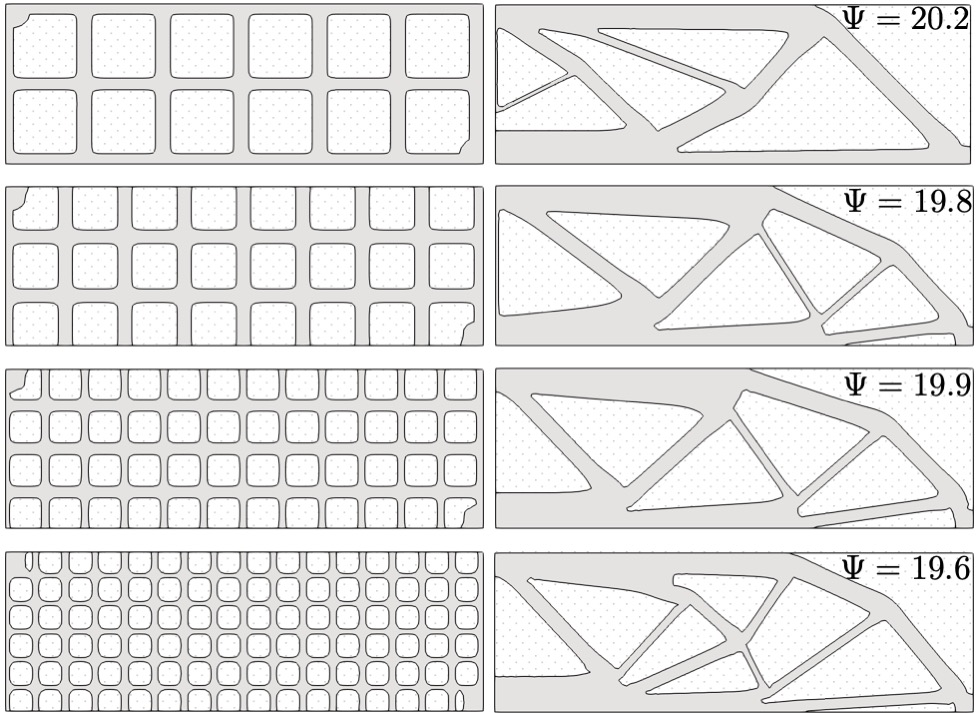}
	\caption{Initial designs with different number of holes seeded (left) and their corresponding final designs (right) for a 2D beam strain energy ($\Psi$) minimization problem.}
	\label{fig:numHolesDependency}
\end{figure}

In some cases, the shortcomings mentioned above can be mitigated by using a ``large enough" number of holes as an initial guess, see \cite{villanueva2014density}. 
Excessive hole seeding, however, requires a fine mesh and a considerable number of design iterations to substantially change the shape of the initial hole pattern and alter conceptually the geometry.
Furthermore, overseeding may produce more refined features without necessarily converging to one specific final design nor impacting performance significantly. 
This case is displayed in Fig. \ref{fig:numHolesDependency} where, despite noticeable differences in the final designs, similar performances are achieved for initial designs using 12, 24, 48, and 96 holes. 

The aforementioned issues have motivated LS approaches capable of nucleating holes. These methods not only circumvent the computational burden present when using a large number of initial holes, but they also eliminate the difficulty of finding a suitable initial seeding pattern. 

Originated from the bubble method (\cite{eschenauer1994bubble}), topological derivatives can be used to seed new holes during the optimization process (\cite{allaire2005structural,burger2004incorporating,wang2004incorporating,norato2007topological,andreasen2019CutFem}). These derivatives evaluate the influence of introducing infinitesimal holes at any point in the design domain during the optimization problem.
Typically, finite size holes are inserted at locations of minimal topological derivative at distinct steps in the optimization process.
Strategies based on topological sensitivity formulations that nucleate holes in regions of low strain energy and stress, among others, have been proposed (\cite{sethian2000structural,belytschko2003topology,sokolowski1999topological,park2008topology}. 
However, regardless of the underlying formulation, using topological derivatives as a guide for hole nucleation introduces user-defined parameters that control the number, frequency, and shape of the seeded holes. 
This can lead to adding arbirtrarily-shaped holes either excessively or scarcely and, as a result, compromise efficiency, robustness and performance (\cite{allaire2004structural}). 
A systematic procedure for deciding if and when nucleating a hole is currently lacking.

Alternative hole seeding mechanisms that use topological derivatives in a non-classical sense, as well as schemes that do not rely on these derivatives, have also been developed.
For example, a topological derivative field can be used to directly extract a geometry if interpreted as a level set field (LSF). 
In \cite{suresh2010199,suresh2013efficient,suresh2013stress,deng2015multi,deng2016multi}, holes are nucleated by redefining the isocontour that separates the material from the void subdomains.
If an implicit LS TO framework is used, the solution of the Hamilton-Jacobi equation can be approximated to allow for the nucleation of holes, as demonstrated by \cite{wang2007extended}.
Alternatively, as shown by \cite{dunning2013new}, a second LSF can be used to construct a pseudo third dimension to nucleate holes in 2D problems.

The class of approaches explored in this paper uses density methods (\cite{bendsoe2004topology}) to nucleate holes in an informed manner.
Since these methods can start from an unbiased design (i.e., homogeneous density distribution) and quickly generate the conceptual layout of an optimal design (\cite{sigmund2013topology}), shortcomings associated with initial hole seeding are avoided.
In these approaches, regions of low density indicate that material is no longer needed and can be used to nucleate holes. Furthermore, sensitivities in the entire design domain characteristic of density-based TO are computed in addition to the shape sensitivities of the LS TO problem, achieving an improved efficiency compared to classical LS TO approaches.

Early attempts to use density methods to advance a well-defined material interface can be found in the context of (concurrent) shape and topology optimization methods. In \cite{kumar1996synthesis}, a ``shape density" field was used to remove parts of the design domain for which densities were under a prescribed threshold. 
\cite{maute1995adaptive,bletzinger1997towards} extensively explored sequentially alternating between shape and topology optimization using separate design and analysis meshes. In their work, post-processed densities define the material interface, which is described by splines (see also \cite{maute1997adaptive,maute1998adaptive}). Smoothing algorithms to realize a crisp material interface from intermediate densities, however, introduces inaccuracies in the description of the geometry. 

More recently, density-based and LS-based TO methods have been combined for purposes other than hole nucleation. In \cite{kang2013integrated}, a classical solid isotropic material with penalization (SIMP) approach is used for material interpolation while a LS approach excludes void domains from the TO process. However, the LSF remains fixed during the design process. 
In an explicit LS TO setting developed by \cite{jansen2019explicit}, feature size control is achieved via geometric constraints on a density field.
The combined LS-density approach in \cite{geiss2018topology,geiss2019combined} couples the LS and density fields to optimize the material distribution within the solid phase. However, the density field is not used for hole nucleation.
 
In this paper, two LS TO approaches that nucleate holes informed by a density field are proposed.
The goal of these two approaches is to seed holes continuously in the optimization process. The location, shape and size of the holes evolve for the specific optimization problem at hand.
In the first approach, a single abstract design variable field governs both the LS and density variables.
In the second approach, two independent design variable fields are introduced to define the LS and density fields, respectively. 
The LSF is one-way coupled with the density field.
Although feature size control based on the density field could be achieved, as demonstrated in \cite{geiss2019level,Barrera2019}, this paper focuses on hole nucleation only.
Note that, similar to the topological derivatives field approach developed by \cite{suresh2010199,suresh2013efficient,suresh2013stress,deng2015multi,deng2016multi}, here an isocontour of a field indicates where a hole needs to be created. 
 However, among other differences, in this paper this indicator density field evolves during the optimization process.

In the material-void problems considered in this stu{\-}dy, a LSF is used to distinguish between phases through crisp, well-defined interfaces. 
This LSF is parametrized by an explicit LS method (\cite{van2013level,sigmund2013topology}).
In addition, the density field interpolates the material properties within the material phase using the SIMP scheme (\cite{bendsoe2004topology}).
The weak form of the governing equations are discretized by the extended finite element method (XFEM, \cite{belytschko2009review}).
Both approaches are studied for linear elastic problems in 2D and 3D, and considering different optimization formulations. 

The remainder of the paper is organized as follows: Section \ref{sec:LSTO} presents the basic components of the LS TO framework employed. Section \ref{sec:HoleNucleationSchemes} explains the two coupling schemes adopted. Section \ref{sec:XFEMFramework} summarizes the discretization/analysis method and describes the general optimization problem formulation. Numerical examples are provided in Section \ref{sec:NumEx}; and Section \ref{sec:Concl} concludes this paper with directions for future work.

\section{Level Set-based Topology Optimization} \label{sec:LSTO}

\subsection{Geometry description} \label{subsec:geomDesc}

The material layout composed of two phases in a design domain, $\Omega_D$, is described by the LSF:
\begin{equation}\label{eq:LSDescription}
\begin{aligned}
\phi(\boldsymbol X)
\begin{cases}
> 0, ~~~~~& \forall~\boldsymbol X \in \Omega_I, \\
< 0, ~~~~~& \forall~\boldsymbol X \in \Omega_{II}, \\
= 0, ~~~~~& \forall~\boldsymbol X \in \Gamma_{I,II},
\end{cases}
\end{aligned}
\end{equation}
where $\Omega_I$ and $\Omega_{II}$ are the material domains of phases $I$ and $II$, respectively, such that $\Omega_D = \Omega_I \cup \Omega_{II}$. The interface, denoted by $\Gamma_{I,II}$, corresponds to the zero LS isocontour, $\phi(\mathbf{X}) = 0$. The problems studied in this work consider solid-void domains where solid is assigned to phase I and void to phase II.

The LSF, bounded between lower, $\phi_{low}$, and upper, $\phi_{up}$, bounds, is a function of a filtered LS design variables field, $\mathcal{\hat{S}^\phi} (\mathbf{X})$, introduced below.

\subsection{LS design variables} \label{subsec:LsDesVars}

In the LS TO framework employed here, the geometry of the solid-void interface is defined by a vector of LS optimization variables, 
${\boldsymbol s}^\phi := \{ {\boldsymbol s}^\phi \in \rm I\!R^{N_s} |~ \phi_{low} \leq {s}^\phi_i \leq \phi_{up}, i=1,...,N_s \}$.
Here, an LS optimization variable is assigned to each node of a structured mesh, and $N_s$ represents the number of nodes in such mesh. 

To increase numerical stability and enhance convergence of the optimization problem, the linear filter presented in \cite{kreissl2012levelset} is applied to ${\boldsymbol{s}}^\phi$ to obtain a vector of filtered LS coefficients, 
$\boldsymbol{\hat{s}}^\phi := \{ \boldsymbol{\hat{s}}^\phi \in \rm I\!R^{N_s} |~\phi_{low} \leq \hat{s}^\phi_i \leq \phi_{up}, i=1,...,N_s \}$.
In this distance-based filter, a filtered LS coefficient, $\hat{s}^\phi_i$, at node $i$, is defined as a function of its neighboring LS optimization variables, ${s}^\phi_j$, at nodes $j$, using the following expression:
\begin{equation}\label{eq:LSlinearFilter}
\begin{aligned}
\hat{s}^\phi_i 
= \frac
{
\displaystyle\sum_{j=1}^{N_{rf}} w_{ij} {s}_j^\phi
}{
\displaystyle\sum_{j=1}^{N_{rf}} w_{ij}
}, 
~
w_{ij} = \max(0,r_f-|\mathbf{X}_i-\mathbf{X}_j|).
\end{aligned}
\end{equation}
The number of nodes within the filter radius, $r_f$ is denoted by $N_{rf}$; and $|\mathbf{X}_i-\mathbf{X}_j|$ is the Euclidean distance between nodes i and j.

The filtered LS coefficients, $\boldsymbol{\hat{s}}^\phi$, are used to parametrize the LS design variable field, $\mathcal{\hat{S}^\phi} ( \mathbf{X})$, following the expression:
\begin{equation}\label{eq:discrLsField}
\begin{aligned}
\mathcal{\hat{S}^\phi} ( \mathbf{X})
= 
\displaystyle\sum_{i=1}^{N^e} \mathcal{N}_i ( \mathbf{X}) ~ \hat{s}_i^\phi,
\end{aligned}
\end{equation}
where
$\mathcal{\hat{S}^\phi} ( \mathbf{X}):=\{ \mathcal{\hat{S}^\phi} ( \mathbf{X}) \in H^1(\Omega_D) |~ \phi_{low} \leq a^\phi \leq \phi_{up}:  a^\phi \in \mathcal{\hat{S}^\phi} ( \mathbf{X})  \}$, 
with the Sobolev space denoted by $H^1$.
The field is interpolated by bi-linear and tri-linear shape functions, $\mathcal{N}( \mathbf{X})$, on quadrilateral and hexahedral meshes in 2D and 3D, respectively. 
The LSF,  $\phi (\mathbf{X})$, is defined as a function of the design variable field, $\mathcal{\hat{S}^\phi} ( \mathbf{X})$, as formulated in Section \ref{sec:HoleNucleationSchemes} separately for each of the hole seeding approaches considered in this paper. 

Unlike implicit approaches where a Hamilton-Jacobi-type equation needs to be solved, the LSF is defined as an explicit function of the LS design variables, ${\boldsymbol s}^\phi$. 
This provides the advantage of solving the optimization problem via a mathematical programming technique. 
The reader is referred to \cite{sigmund2013topology} and \cite{van2013level} for more details.

\subsection{LS regularization} \label{subsec:LsReg}
To avoid spurious oscillations in the LSF, the regularization scheme of \cite{geiss2019regularization} is adopted. This approach promotes a uniform spatial gradient of the LSF near the solid-void interface while converging to either a positive or negative target value away from the interface. 
A truncated signed distance field used as target field, $\tilde\phi$, is constructed via the heat method (\cite{crane2017heat}). This entails solving two linear finite element problems (i.e., a heat conduction problem and a Poisson's problem) on the entire design domain.
The target field is enforced through a penalty in the objective, $P_{Reg}$, of the following form:
\begin{equation}\label{eq:LsRegEq}
\begin{aligned}
	P_{Reg} = 
	w_{\phi} \frac{\displaystyle\int_{\Omega_D} (\phi-\tilde\phi)^2 dV}{\displaystyle\int_{\Omega_D} \tilde\phi_{Bnd}^2 dV} + 
	w_{\nabla\phi} \frac{\displaystyle\int_{\Omega_D} | \nabla \phi- \nabla\tilde\phi|^2 dV}{\displaystyle\int_{\Omega_D} dV},
\end{aligned}
\end{equation}
where $\tilde\phi_{Bnd}$ is the difference between the upper, $\tilde\phi_{up}$, and lower, $\tilde\phi_{low}$, bounds in the target LSF. 
In \cite{geiss2019regularization}, the weights $w_{\phi}$ and $w_{\nabla\phi}$ were kept constant in the entire design domain. In the current work, to balance the influence of the regularization components in the vicinity and away from the material interface, these weights are defined as:
\begin{equation}\label{eq:LsRegWeights}
\begin{aligned}
&	w_{\phi} = 
	w_{\phi_1} \alpha + w_{\phi_2} (1-\alpha),  \\
&	w_{\nabla\phi} = 
	w_{\nabla\phi_1} \alpha + w_{\nabla\phi_2}(1-\alpha).
\end{aligned}
\end{equation}
The weights $w_{\phi_1},w_{\phi_2}$ control the mismatch between design and target LSF near and far from the material interface, respectively. The weights $w_{\nabla\phi_1},w_{\nabla\phi_2}$ do the same for the LS spatial gradient component in Eq. \ref{eq:LsRegEq}. The parameter $\alpha$, defined as:
\begin{equation}\label{eq:alphaLsReg}
\begin{aligned}
	\alpha(\tilde\phi) &= 
	e^{-\gamma_{P_{Reg}} (\tilde\phi/\phi_{Bnd})^2},
\end{aligned}
\end{equation}
controls the region of influence of the regularization through $\gamma_{\phi Reg}$, as seen in Fig. \ref{fig:LsRegAlphaParam}. 
Larger $\gamma_{\phi Reg}$ values decrease $\alpha$ and thus, the weight of the regularization in the vicinity of the interface.

\begin{figure}[h!]
	\center
	\includegraphics[width=0.5\linewidth]{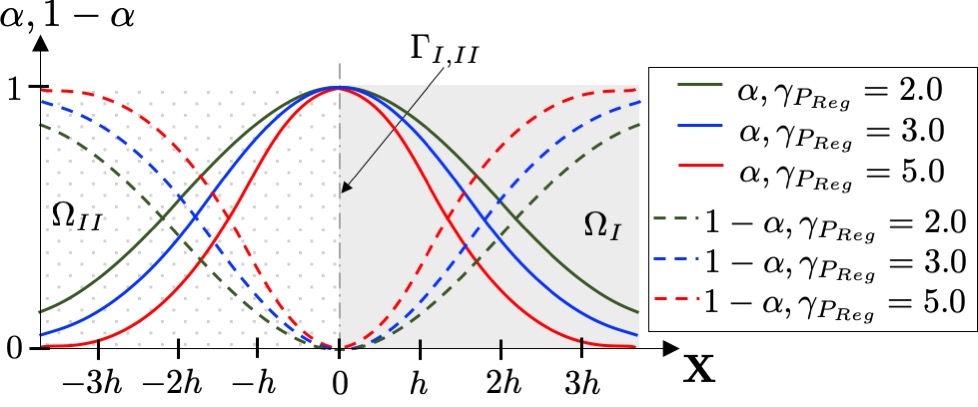}
	\caption{LS regularization penalty as function of the $\alpha$ parameter in the vicinity of the interface $\Gamma_{I,II}$ for a mesh with element length $h$.}
	\label{fig:LsRegAlphaParam}
\end{figure}
%

\section{Hole Nucleation via Density Methods} \label{sec:HoleNucleationSchemes}

Density methods start typically from an unbiased, i.e., homogeneous, density distribution and quickly converge to the conceptual layout of the optimized design. This attractive feature is used here to instruct the LS optimization problem to remove material where it is not needed.
To this end, we introduce a nodal fictitious density field, $\rho(\mathbf{X})$ bounded between 0 and 1 in the LS TO problem to nucleate holes in a mechanically-informed manner. 

Following the same scheme presented in the previous section for the LSF,  $\rho(\mathbf{X})$ is defined as an explicit function of a filtered density design variables field, 
$\mathcal{\hat{S}^\rho} ( \mathbf{X}):=\{ \mathcal{\hat{S}^\rho} ( \mathbf{X}) \in H^1(\Omega_D) |~ 0 \leq a^\rho \leq 1:  a^\rho \in \mathcal{\hat{S}^\rho} ( \mathbf{X})  \}$; 
i.e.,  $\rho(\mathbf{X}):= \rho ( \mathcal{\hat{S}}^\rho(\mathbf{X}) )$. 
The field $\mathcal{\hat{S}}^\rho(\mathbf{X})$ is discretized using Eq. \ref{eq:discrLsField} and a vector of filtered density coefficients, 
${\boldsymbol{\hat{s}}}^\rho := \{ {\boldsymbol{\hat{s}}}^\rho \in \rm I\!R^{N_s} |~ 0 \leq \hat{s}^\rho_i \leq 1, i=1,...,{N_s} \}$. 
Similarly, the linear filter in Eq. \ref{eq:LSlinearFilter} is employed to define $\boldsymbol{\hat{s}}^\rho$ in terms of a vector of nodal density optimization variables, 
${\boldsymbol s}^\rho := \{ {\boldsymbol s}^\rho \in \rm I\!R^{N_s} |~ 0 \leq {s}^\rho_i \leq 1, i=1,...,{N_s} \}$. 

In this section we present two approaches to couple the LS and density design variables. 
In both cases, the couplings of the LS and density fields are constructed to achieve that (i) a positive LSF is associated with high densities to represent material, i.e., $\rho(\mathbf{X}) \rightarrow 1 \implies \phi(\mathbf{X}) > 0$; and (ii) a negative LSF is linked to densities close to zero to identify void, i.e., $\rho(\mathbf{X}) \rightarrow 0 \implies \phi(\mathbf{X}) < 0$. 
In the first approach, both the LS and density fields, $\phi(\mathbf{X})$ and $\rho(\mathbf{X})$, are both functions of a single set of optimization variables; i.e., ${\boldsymbol s}^\phi={\boldsymbol s}^\rho$. We identify this approach as Single-Field Coupling (SFC).
The second approach interpolates $\phi(\mathbf{X})$ and $\rho(\mathbf{X})$ using two independent sets of optimization variables. This approach is termed Two-Field Coupling (TFC). 

\subsection{Single-field coupling strategy} \label{subsec:SingleFieldCoupling}	

In the first proposed, we assume a single vector of ${N_s}$ abstract design variables, $\boldsymbol{s}=\boldsymbol{s^\phi}=\boldsymbol{s^\rho}$, to compute the vectors of filtered coefficients, $\boldsymbol{\hat{s}}=\boldsymbol{\hat{s}^\phi}=\boldsymbol{\hat{s}^\rho}$, using Eq. \ref{eq:LSlinearFilter}, and interpolate $\mathcal{\hat{S}}(\mathbf{X}) =\mathcal{\hat{S}^\phi}(\mathbf{X}) = \mathcal{\hat{S}^\rho}(\mathbf{X})$ using Eq. \ref{eq:discrLsField}. 
Here, the LSF is defined by the following linear relation:
\begin{equation}\label{eq:SFC_AdvToLsVars}
\begin{aligned}
	\phi(\mathbf{X})= \phi_{rt} (\phi_{sh} - \mathcal{\hat{S}}(\mathbf{X})).
\end{aligned}
\end{equation}
Similarly, the density field is obtained by using:
\begin{equation}\label{eq:SFC_AdvToDensVars}
\begin{aligned}
	\rho(\mathbf{X}) =
	\begin{cases}
	\displaystyle\frac{\mathcal{\hat{S}}(\mathbf{X}) -\phi_{sh}}{1-\phi_{sh}}, & \forall~\mathbf{X}:\phi(\mathbf{X})\geq 0, \\
	not~defined, & \forall~\mathbf{X}:\phi(\mathbf{X})<0.
	\end{cases}
\end{aligned}
\end{equation}
In both Eqs. \ref{eq:SFC_AdvToLsVars} and \ref{eq:SFC_AdvToDensVars},  $\phi_{sh}$ is a threshold that alters the LS value used to differentiate between solid and void. In this work, $\phi_{sh}$ is kept constant, i.e., $\phi_{sh} = 0.5$. 
The scaling parameter $\phi_{rt}$ in Eq. \ref{eq:SFC_AdvToLsVars} takes into account the mesh size, and is set to $\phi_{rt}\approx [3h, 5h]$, where $h$ is the element length.  
The density field described by Eq. \ref{eq:SFC_AdvToDensVars} is only defined in the material subdomain, $\Omega_{I}$.

\subsection{Two-field coupling strategy} \label{subsec:TwoFieldCoupling}

The second coupling strategy uses two sets of independent LS and density optimization variables, $\boldsymbol{s^\phi}$ and $\boldsymbol{s^\rho}$.
Consequently, the optimization process in this case operates on $2N_s$ optimization variables; i.e., 
$\boldsymbol{s}=[\boldsymbol{s^\phi}, \boldsymbol{s^\rho}]$.
Here, the LS and density fields are assumed to be identical to the LS and density design variable fields, respectively; i.e., $\phi(\mathbf{X})= \mathcal{\hat{S}}^\phi(\mathbf{X})$ and $\rho(\mathbf{X}) = \mathcal{\hat{S}}^{\rho}(\mathbf{X})$.

In this scheme, holes are seeded in the course of the optimization process by promoting a negative LSF in regions of low density through a penalty term, $\bar{p}_{\rho\phi}(\mathbf{X})$, defined as:
\begin{equation}\label{eq:TFCNonSmoothForm}
\begin{aligned}
	\bar{p}_{\rho\phi}(\mathbf{X}) &=
	\begin{cases}
	0, & 			\forall~\rho\geq\rho_{th} \\
	\max \left( 0, \displaystyle\frac{\phi -\phi_{th}}{\phi_{up}-\phi_{th} } \right), & \forall~\rho<\rho_{th}.
	\end{cases}
\end{aligned}
\end{equation}
This coupling penalty is active in a region bounded by the $\phi_{th}$ and $\rho_{th}$ thresholds, described below. 
High LS values in regions of low densities are penalized in a decreasing manner from 1 at $\phi=\phi_{up}$ to 0 for $\phi \leq \phi_{th}$.
To serve as a hole nucleation mechanism, it is intentionally formulated by a discontinuous function that has zero gradients with respect to the density optimization variables and exhibits non-zero gradients with respect to the LS optimization variables only below $\rho_{th}$. 

\subsubsection{LS and density thresholds} \label{subsubsec:TFCThs}

The LS threshold, $\phi_{th}$, deactivates the penalization for $\phi(\mathbf{X})<\phi_{th}$. 
This threshold is set to a value below zero, i.e., $ [0.10,0.25] \phi_{low}$. 
Setting $\phi_{th} \geq 0$ would result in issues associated with the robustness of the hole nucleation process.
An explanation of these issues is provided in Section \ref{subsubsec:holeNucProc}.

The $\rho_{th}$ threshold is a parameter that decreases from an initial value, $\rho^0_{th}$, to zero during the optimization process.
A decreasing $\rho_{th}$ is desired since intermediate densities at a later stage of the design process could promote the creation of spurious holes, see Section \ref{subsubsec:densShift}.
Note that $\rho_{th}^0$ should be lower than the prescribed initial density, $\rho_0$, in the first optimization iteration.
Setting $\rho^0_{th}>\rho_0$ would activate the penalty in the entire domain and could result in moving the entire LSF from material to void.

Through a continuation scheme active in the first $\mathcal{D}_{c}$ design iterations, the density threshold, $\rho_{th}$, is computed at every continuation step using:
\begin{equation}\label{eq:expFuncContDec}
\begin{aligned}
	\rho_{th} = \rho_{th}^0
	\begin{cases}
	\displaystyle 1 - \left( \frac{\mathcal{D}_{it}}{\mathcal{D}_{c}}\right)^{\eta_{\rho_{th}}}, & \forall~\mathcal{D}_{it} \leq \mathcal{D}_{c} \\
	0, & \forall~\mathcal{D}_{it} > \mathcal{D}_{c},
	\end{cases}	
\end{aligned}
\end{equation}
where $\mathcal{D}_{it}$ is the design iteration index.
As seen in Fig. \ref{fig:DensThresholdScheme}, $\rho_{th}$ decreases to gradually inhibit the nucleation of holes.
The coupling is inactive for $\mathcal{D}_{c}<\mathcal{D}_{it}$.
The continuation step size, $\mathcal{D}_{st}$, is assumed to be constant.
The exponent $\eta_{\rho_{sh}}$ controls the rate at which the shift increases in each continuation step, and is set to $\eta_{\rho_{th}}\approx [1.5-2.0]$ to attenuate changes in this threshold at early stages of the continuation scheme. 

\begin{figure}[h!]
	\center
	\includegraphics[width=0.5\linewidth]{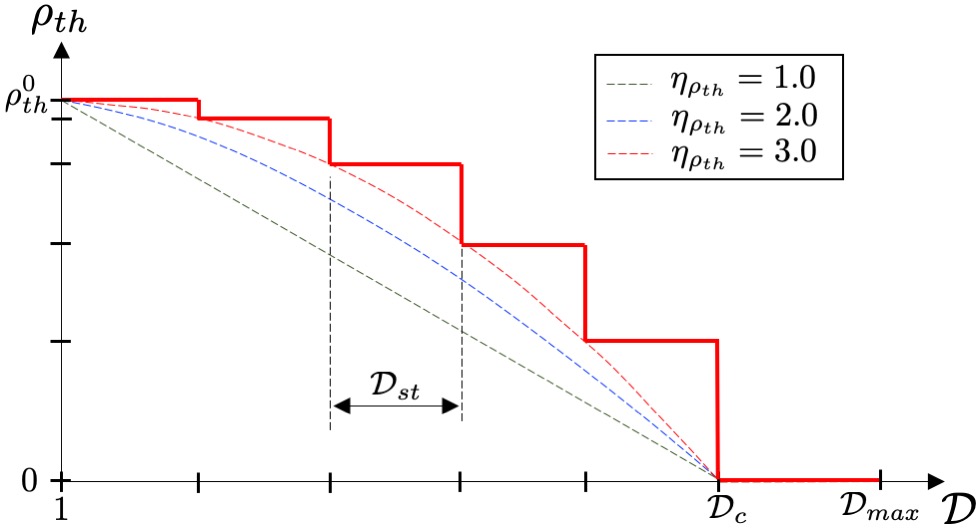}
	\caption{Continuation scheme in Eq. \ref{eq:expFuncContDec} for the density threshold of the TFC approach.}
	\label{fig:DensThresholdScheme} 
\end{figure}

\subsubsection{Smooth coupling formulation}

\begin{figure}[h!]
	\center
	\includegraphics[width=0.5\linewidth]{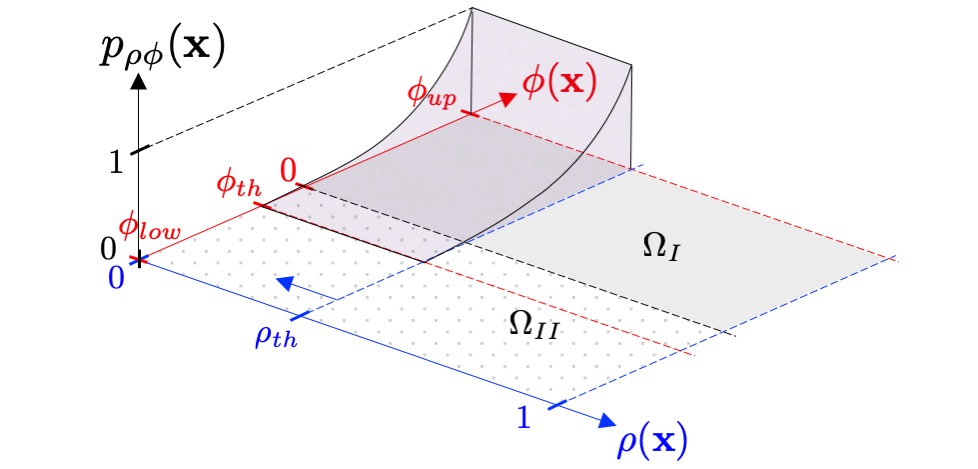}
	\caption{Landscape of the smooth coupling penalization in Eq. \ref{eq:TFCFomAllDomain} for the TFC approach.}
	\label{fig:TFCScheme} 
\end{figure}

To obtain a differentiable penalty formulation with respect to the LSF, the expression in Eq. \ref{eq:TFCNonSmoothForm} is approximated by a smooth function as follows:
\begin{equation}\label{eq:TFCFomAllDomain}
\begin{aligned}
	{p}_{\rho\phi}(\mathbf{X}) &=
	\begin{cases}
	0, ~
	\forall~\rho\geq\rho_{th} 
	\\
	\displaystyle\frac{
	\left\{ \left[ \max \left( 0, \displaystyle\frac{\phi -\phi_{th}}{\phi_{up}-\phi_{th} } \right) \right]^2 
	+ \xi^2 \right\}^{\frac{1}{2}} - \xi
	}{
	\left( 1 + \xi^2 \right)^{\frac{1}{2}}  - \xi
	}, & 			
	\\ 
	\forall~\rho<\rho_{th}.
	\end{cases} \\
\end{aligned}
\end{equation}
The nondimensional parameter $\xi$ smoothes the transition of the penalty, and is typically between $\xi\approx [0.1-1.0]$. 
An illustration of this smooth two-field density-LS coupling scheme is presented in Fig. \ref{fig:TFCScheme}. 

The TFC strategy introduced here represents an improvement over the definition of the coupling penalty in \cite{geiss2019level}. Numerical studies showed that the approach in \cite{geiss2019level} is prone to oscillations of the LSF in low density regions. The spurious interaction between density and LS fields is avoided by the discontinuous formulation of the penalty term in Eq. \ref{eq:TFCFomAllDomain}. 
Since the gradients of the penalty term with respect to the density optimization variables are zero, this formulation does not provide information that would promote the removal of a hole through increasing the local densities. 
 However, a hole can still be removed through shape changes. Thus, the coupling penalty formulation in Eq. \ref{eq:TFCFomAllDomain} should be considered a hole nucleation scheme.

Note that the choice of variables in this coupling strategy is not unique. Either the sets of LS and density optimization variables, or any choice of filtered and/or projected quantities derived from them could be discretized and coupled. Our choice of using $\mathcal{\hat{S}}^{\phi}(\mathbf{X})$ and $\mathcal{\hat{S}}^{\rho}(\mathbf{X})$ is motivated by an increased in resolution and design freedom (compared to elemental variables), as well as smoothness in the optimization problem. 

\subsubsection{Hole nucleation process} \label{subsubsec:holeNucProc}

\begin{figure}[h!]
	\center
	\includegraphics[width=0.5\linewidth]{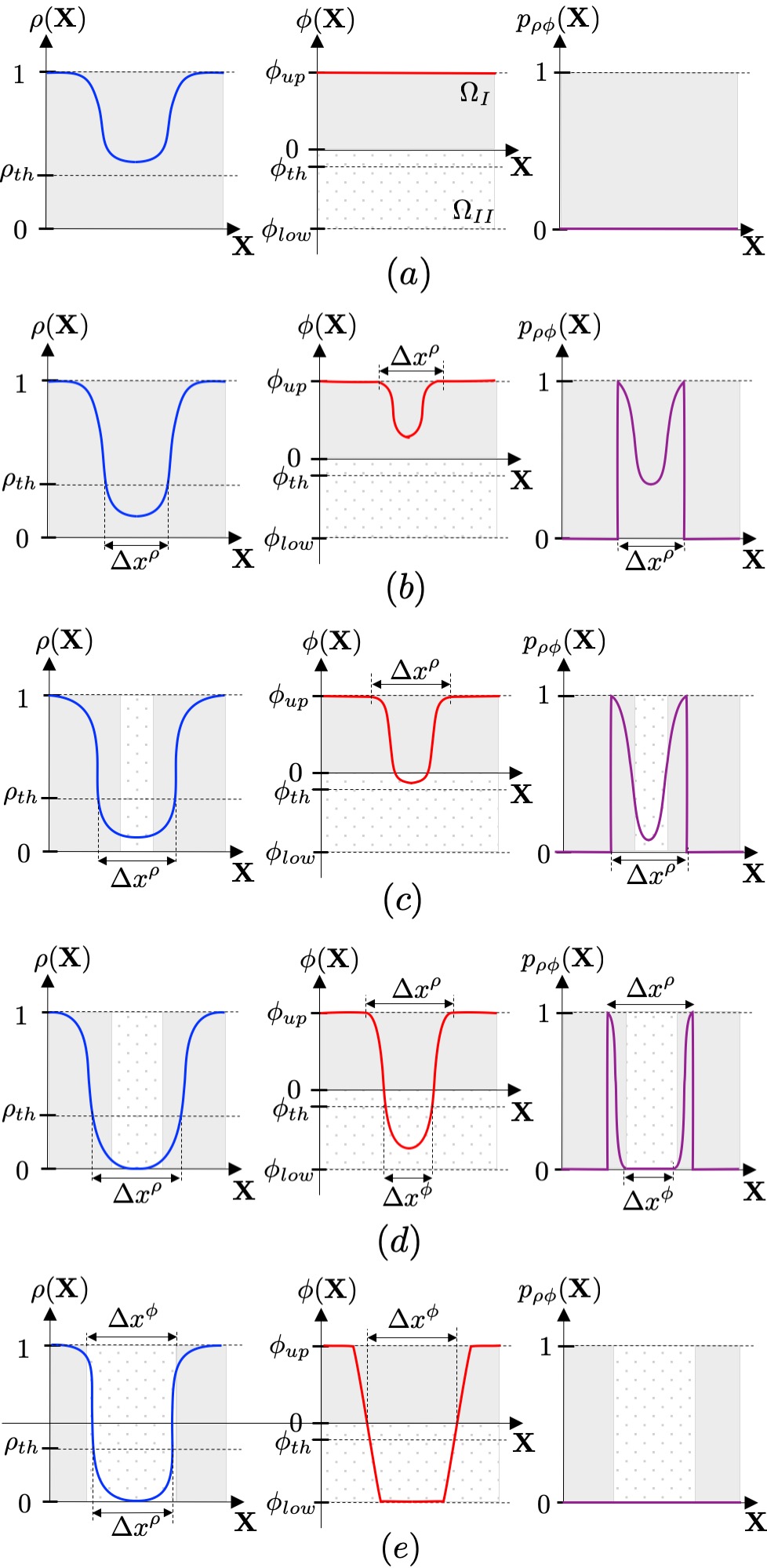}
	\caption{Evolution of the density, LS, and penalty fields at all stages of the nucleation process via the TFC approach.}
	\label{fig:TFC1DExample}
\end{figure}

Figure \ref{fig:TFC1DExample} depicts the effect of the TFC penalty in a simple 1D example. 
The density, LS, and coupling penalty fields, i.e., $\rho(\mathbf{X})$, $\phi(\mathbf{X})$, and ${p}_{\rho\phi}(\mathbf{X})$, are plotted over $\mathbf{X}$ at different stages of the nucleation process of a hole. 
Initially, in Fig. \ref{fig:TFC1DExample}(a), the coupling penalty is inactive in the material domain, i.e., $\phi(\mathbf{X})>0$ and $\rho(\mathbf{X})>\rho_{th}$ everywhere. 
As the local density value decreases, the coupling penalty is activated along $\Delta x^{\rho}$ defined by the region of the density field smaller than $\rho_{th}$. As a result, the LSF in $\Delta x^{\rho}$ is lowered as seen in Fig. \ref{fig:TFC1DExample}(b). 
The nonlinear distribution of the penalty (see Eq. \ref{eq:TFCFomAllDomain}) is maximum at $\phi(\mathbf{X})=\phi_{up}$, and its effect diminishes as the LSF approaches $\phi_{th}$. 
Eventually, in Fig. \ref{fig:TFC1DExample}(c), a hole is nucleated by the coupling penalty, which continues decreasing the LSF (to a lesser extent) while it is above $\phi_{th}$.
Note that, as mentioned in Section \ref{subsubsec:TFCThs}, using $\phi_{th} \geq 0$ can compromise the hole nucleation process. 
Setting $\phi_{th}=0$ may introduce a robustness issue since the LSF could oscillate between positive and negative values, nucleating and removing small holes intermittently.
If $\phi_{th}>0$, the coupling terminates prematurely and the LSF is pushed back into the material domain before crossing the zero isocontour, preventing a hole from being nucleated.
This is a consequence of the LS regularization scheme detailed in Section \ref{subsec:LsReg}.
 Hence, $\phi_{th}<0$ sufficiently away from the zero isocontour is used in this work.
In Fig. \ref{fig:TFC1DExample}(d), the LSF below $\phi_{th}$ along $\Delta x^{\phi}$ is no longer affected by the penalty, but continues decreasing until it reaches $\phi_{low}$ due to the effect of the LS regularization. At this stage, the density field is converged and the penalty is acting only in a small portion of the domain contained under $\rho_{th}$ and above $\phi_{th}$ (i.e., $p_{\rho\phi}(\mathbf{X})>0$ in $\Delta x^{\rho} \cap \Delta x^{\phi}$). 
Finally, Fig. \ref{fig:TFC1DExample}(e) shows the regularized LSF with a hole nucleated in a region of low density. The penalization is no longer active since the process has been completed.
A study of the influence of the two thresholds, $\rho_{th}$ and $\phi_{th}$ is presented in Section \ref{subsubsec:SupportStruc2DObjAndSensTFC}.

%

\subsection{Single versus two-field coupling} \label{subsec:SFCvTFC}

The two coupling strategies presented above enable hole seeding during the optimization process. 
In both cases, the design process consists of a handover from a pure density problem that does not need hole seeding but can produce intermediate density at convergence to a pure LS problem that relies on hole nucleation but features a crisply defined material interface.
In terms of the gradients, only design sensitivities with respect to the density variables exist everywhere within the solid domain until holes are nucleated. Once that happens, non-zero LS shape sensitivities are also present in the vicinity of the material interface.
Non-zero density sensitivities in an LS TO approach not only accelerate convergence, but also facilitate convergence to (local) minima with satisfactory performance and mitigate issues associated with initial seeding (see Fig. \ref{fig:holesArrangementDependency}). 
Note that non-zero sensitivities in the entire design domain can also exist as a result of the LS regularization method used. The method adopted here is one of such cases (see Section \ref{subsec:LsReg}). 

In both hole nucleation schemes, the computational costs associated with physical analyses are reduced. The mesh resolution is no longer dominated by the need to resolve an initial hole seeding. Hence, the proposed method might allow using coarser meshes depending on the feature sizes in the optimized design. In the early stages of the optimization process, complex intersection configurations due to user-defined initial hole seeding are avoided. This reduces the number of intersected elements and thus, the computational cost when using immersed boundary techniques for the physical analyses.
Furthermore, the number of total design iterations might decrease due to faster convergence.
However, in contrast to evaluating LS shape sensitivities by only considering intersected elements, the contributions of all elements to the sensitivity equations within the solid domain need to be computed.
The increased cost of the TFC approach over the SFC formulation associated with the increased number of optimization variables is insignificant because of the class of optimization algorithms and the type of sensitivity analysis used in this manuscript; see Section \ref{sec:NumEx}.

\subsection{Material interpolation via augmented SIMP} \label{subsec:DensTO}

Transition regions between high and low density material are commonly found in density-based methods. For example, SIMP-like approaches typically require a suitable formulation of the optimization problem and large penalties to converge to [0-1] material distributions. However, a penalization effect on intermediate densities, and thus their complete removal from the design domain, cannot be guaranteed (\cite{sigmund2013topology}).

\subsubsection{Shifted density field} \label{subsubsec:densShift}

To avoid intermediate densities, we introduce a density shift, $\rho_{sh}$, to compute a shifted density field, $\tilde\rho(\mathbf{X})$, using the following expression:
\begin{equation}\label{eq:densShift}
\begin{aligned}
\tilde\rho(\mathbf{X})
=
\rho_{sh}
+
(1-\rho_{sh}) \rho(\mathbf{X}).
\end{aligned}
\end{equation}
%
The parameter $\rho_{sh}$ increases from an initial value, $\rho^0_{sh}$, to 1 during the optimization process.
It has the effect of shifting the minimum density in the material domain, as shown in Fig. \ref{fig:SFCScheme} for the SFC approach. The larger $\rho_{sh}$, the more restricted the density distribution becomes. If $\rho_{sh}=1$, a uniform density field is obtained.
\begin{figure}[h!]
	\center
	\includegraphics[width=0.5\linewidth]{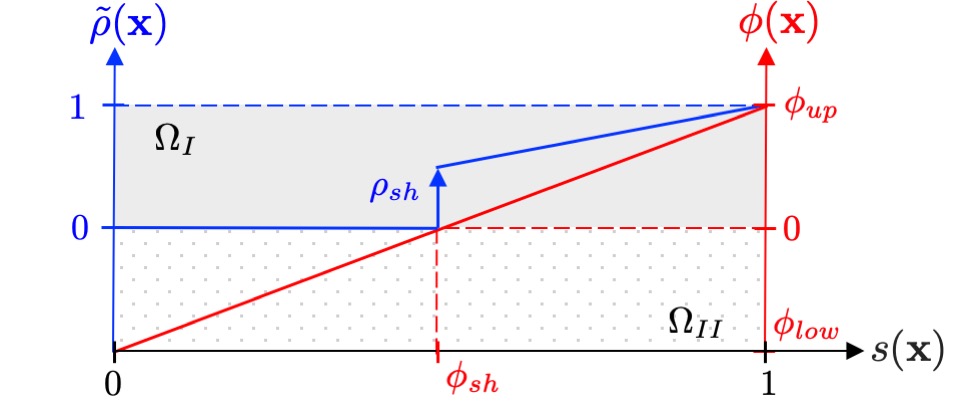}
	\caption{Effect of density shift, $\rho_{sh}$, in the SFC approach.}
	\label{fig:SFCScheme}
\end{figure}

By shifting the densities, not only a constant density field is achieved in the final design, but ill-conditioning of the discretized governing equations due to large differences in material properties is alleviated. 
The latter is achieved because areas/volumes of low density and low material stiffness are removed by controlling the minimum density in the solid domain.
To attain well-conditioned systems, $\rho^0_{sh} \approx [0.1-0.2]$ unless a smaller value is required to satisfy a mass or volume constraint. 
Moreover, the density shift boosts the LS sensitivities of new holes since, by
advancing the shift, the problem transitions from a pure SIMP to a pure LS formulation. Thus, the density sensitivities decrease, and the LS shape sensitivities increase.

A continuation scheme similar to the one used for the density threshold (see Eq. \ref{eq:expFuncContDec}), is employed to update the density shift:
\begin{equation}\label{eq:expFuncContInc}
\begin{aligned}
	\rho_{sh} = 
	\begin{cases}
	\displaystyle \rho_{sh}^0 + ( 1 - \rho_{sh}^0 ) \left( \frac{\mathcal{D}_{it}}{\mathcal{D}_{c}}\right)^{\eta_{\rho_{sh}}}, & \forall~\mathcal{D}_{it} \leq \mathcal{D}_{c} \\
	1, & \forall~\mathcal{D}_{it} > \mathcal{D}_{c}.
	\end{cases}
\end{aligned}
\end{equation}
Eq. \ref{eq:expFuncContInc} ensures the elimination of intermediate densities for $\mathcal{D}_{it} > \mathcal{D}_{c}$, see Fig. \ref{fig:DensShiftScheme}.
Every $\mathcal{D}_{st}$ iterations, the density shift is updated from an initial value, $\rho^0_{sh}$, to 1.
Numerical studies suggest that setting the exponent $\eta_{\rho_{sh}}\approx [1.5-2.0]$, together with an appropriate number of continuations steps (i.e., $\mathcal{D}_{st}/\mathcal{D}_{c}\approx[5-10]$) sufficiently mitigate the effect of jumps in the densities and promote a smooth evolution of the optimization problem. 

\begin{figure}[h!]
	\center
	\includegraphics[width=0.5\linewidth]{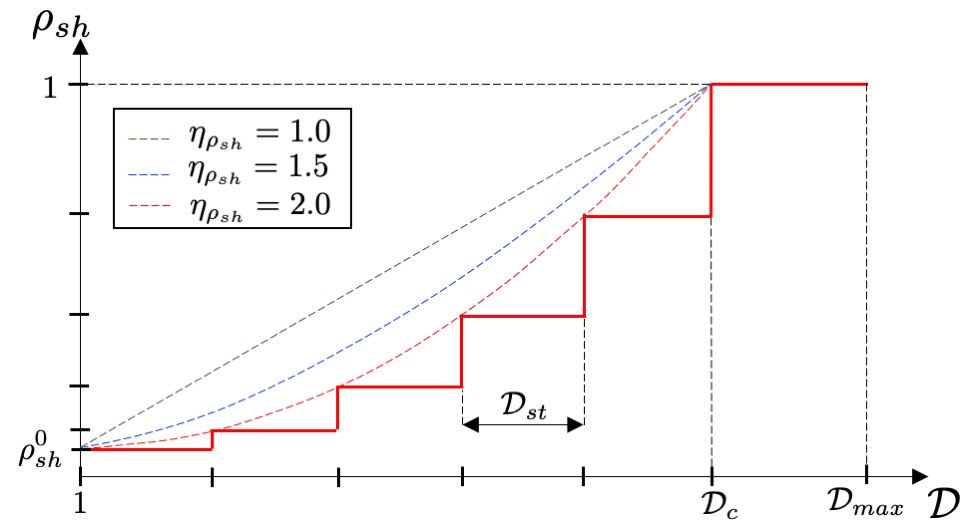}
	\caption{Continuation scheme for density shift, $\rho_{sh}$.}
	\label{fig:DensShiftScheme}
\end{figure}

\subsubsection{Material interpolation} \label{subsubsec:matInterp}

The classical SIMP approach, as presented by \cite{bendsoe2004topology}, is employed to relate the shifted fictitious density field, $\tilde\rho(\mathbf{X})$, to physical material properties. However, alternatives such as the Rational Approximation of Material Properties (RAMP) method (\cite{stolpe2001alternative}) are equally suitable.
The material density, $\theta(\mathbf{X})$, is interpolated using:
\begin{equation}\label{eq:varyingMatDensForm}
\begin{aligned}
	\theta(\mathbf{X}) = \theta_0  ~ \tilde\rho(\mathbf{X});	
\end{aligned}
\end{equation}
and the Young's modulus, $E(\mathbf{X})$, is computed using the following power law:
\begin{equation}\label{eq:varyingMatPropForm}
\begin{aligned}
	E(\mathbf{X}) = E_0 \left(\tilde\rho(\mathbf{X}) \right)^{\beta_{\rho}}.	
\end{aligned}
\end{equation}
The properties of the bulk material density and Young's modulus are denoted by $\theta_0$ and $E_0$. The Young's modulus is penalized by the SIMP exponent, denoted by $\beta_{\rho}$.
A constant density within an element is interpolated by averaging the shifted nodal densities in such element.

While material properties are interpolated using a shifted filtered density field, $\tilde\rho(\mathbf{X})$, the hole seeding strategies presented earlier in this section couple the unshifted filtered density field, $\rho(\mathbf{X})$, to the LSF, $\phi(\mathbf{X})$. 
In the TFC approach, using the shifted densities in the formulation of the penalty term would terminate the coupling too early since hole nucleation capabilities are lost once the density shift exceeds the density threshold, i.e., $\rho_{sh} > \rho_{th} \implies p_{\rho\phi}(\mathbf{X})=0$.
Also, later in the optimization process, $\tilde\rho(\mathbf{X}) \to 1$ does not imply $\rho(\mathbf{X}) \to 1$. Hence, intermediate densities might still be present in $\rho(\mathbf{X})$, without affecting the interpolation of the material properties and thus, the optimization problem. 
For this reason, the density threshold is decreased gradually approaching zero (see Section \ref{subsubsec:TFCThs}).
Otherwise, spurious holes might be nucleated in regions of high shifted densities in designs that are almost converged.
A motivation for using the augmented material interpolation scheme presented in this section, together with an analysis of its influence, is discussed in Section \ref{subsubsec:SupportStruc2DIntermediateDens}.

\section{Optimization Framework} \label{sec:XFEMFramework}


The physical responses of the systems are predicted by the XFEM in this paper. 
Since LS TO is not restricted to a particular immersed boundary technique, alternative approaches such as CutFEM (\cite{burman2015cutfem,burman2019cut}) or HIFEM (\cite{soghrati2016application}) could also be used for physical and sensitivity analyses. 

A generalized Heaviside enrichment strategy (\cite{makhija2014numerical,terada2003finite,tran2011multiple}) is employed to avoid artificial coupling of disconnected material. The response is consistently interpolated in elemental subdomains with the same phase.
The discretized state variable field, $\mathbf{\hat{u}}_i(\mathbf{X})$, at node $i$ of a two material problem ($\Omega_I$, $\Omega_{II}$) is approximated as:
\begin{equation}\label{eq:HeavisideEnrich}
\begin{aligned}
	\mathbf{\hat{u}}_i(\mathbf{X})
	=
	\sum^{L}_{l=1}
	\bigg(
	H(-\phi(\mathbf{X}))\sum^{N_N^e}_{k=1} \mathcal{N}_k(\mathbf{X})\delta_{lq}^{k}\mathbf{{u}}_{il}^{k,\Omega_{I}} \\
	+ 
	H(\phi(\mathbf{X}))\sum^{N_N^e}_{k=1} \mathcal{N}_k(\mathbf{X})\delta_{lq}^{k}\mathbf{{u}}_{il}^{k,\Omega_{II}}
	\bigg),
\end{aligned}
\end{equation} 
where the Heaviside function, $H$, is determined by the LSF as:
\begin{equation}\label{eq:HeavisideFunc}
\begin{aligned}
	H(\phi) =
	\begin{cases}
		1, & \forall ~\phi(\mathbf{X}) > 0 \\
		0, & \forall ~\phi(\mathbf{X}) < 0.
	\end{cases}
\end{aligned}
\end{equation}
The maximum number of enrichment levels is denoted by $L$, $\mathcal{N}_k(\mathbf{X})$ is the nodal shape function and $\delta_{lq}^{k}$ is the Kronecker delta which selects the active enrichment level $q$ for node $k$. $\delta_{lq}^{k}$ ensures that displacements at node $k$ are only interpolated by a single set of degrees of freedom (DOFs), $\mathbf{{u}}_{il}^{k}$, such that the partition of unity principle is satisfied (\cite{babuvska1997partition}). The number of nodes per element is denoted by $N_N^e$. For more details about the generalized Heaviside enrichment strategy used in this work, the reader is referred to \cite{makhija2014numerical}. The described XFEM framework has been successfully applied to various TO multiphysics problems (e.g., in \cite{maute2015level,coffin2016level,villanueva2017cutfem,behrou2017level,pizzolato2017design}). 

\subsection{Governing equations} \label{subsec:GovEqs}

The following augmented residual equation with stabilization terms is considered in this work:
\begin{equation}\label{eq:resStrucGovEq}
\begin{aligned}
	\mathbf{R}
	=
	\mathbf{R}^U + 
	\mathbf{R}_\Gamma^N  +
	\mathbf{R}_\Gamma^G  +
	\mathbf{R}^S,
\end{aligned}
\end{equation}
in which the weak form of the linear elastic governing equation, $\mathbf{R}^{U}$, is defined as:
\begin{equation}\label{eq:ResEq}
\begin{aligned}
	\mathbf{R}^U = 	
	\int_{\Omega_{I}} \delta \boldsymbol{\epsilon} : \boldsymbol{\sigma}~ dV - 
	\int_{\Gamma_N^{\Omega_I}} \delta \mathbf{u}~\mathbf{T}_N dA.
\end{aligned}
\end{equation}
The displacement field and the test function are denoted by $\mathbf{u}$ and $\delta\mathbf{u}$, respectively. Traction forces, $\mathbf{T}_N$, are applied along the boundary $\Gamma_N^{\Omega_I}$. 
The material tensor for isotropic linear elasticity, $\mathbf{D}$, together with the infinitesimal strain tensor, $\boldsymbol{\epsilon}=\frac{1}{2}(\nabla\mathbf{u} + \nabla\mathbf{u}^T)$, define the Cauchy stress $\boldsymbol{\sigma} = \mathbf{D}: \boldsymbol{\epsilon}$.

The $\mathbf{R}_\Gamma^N$, $\mathbf{R}_\Gamma^G$, and $\mathbf{R}^S$ terms correspond to the weakly enforced essential boundary and interface conditions, face-oriented ghost stabilization, and selective structural springs, respectively.

\subsubsection{Weak enforcement of boundary and interface conditions}

\label{sec:Nitsche}
The unsymmetric version of Nitsche's method (\cite{nitsche1971variationsprinzip}) is employed to weakly enforce Dirichlet boundary and interface conditions (\cite{burman2012fictitious}). These conditions are applied using:
\begin{equation}\label{eq:NitscheFrom}
\begin{aligned}
	\mathbf{R}_\Gamma^N 
	= &
	 - \int_{\Gamma} 
	\llbracket \delta \mathbf{u} \rrbracket ~ \boldsymbol\sigma \cdot \mathbf{N} dA 
	+ \int_{\Gamma} 
	\delta ( \boldsymbol\sigma \cdot \mathbf{N} ) \llbracket \mathbf{u} \rrbracket dA 
	\\
	& + \gamma_N~E / h\int_{\Gamma} 
	\llbracket \delta \mathbf{u} \rrbracket \llbracket \mathbf{u} \rrbracket dA.
\end{aligned}
\end{equation}
The outward normal vector on a domain or interface boundary is denoted by $\mathbf{N}$.
The jump operator $\llbracket \boldsymbol{\cdot} \rrbracket$ with respect to a target state variable $\bar{\mathbf{u}}$ is defined as:
\begin{equation}\label{eq:jumpOp}
\begin{aligned}
	\llbracket \mathbf{u} \rrbracket = \mathbf{u} - \bar{\mathbf{u}}, ~~~~~
	\llbracket \delta \mathbf{u} \rrbracket = \delta \mathbf{u} - \delta \bar{\mathbf{u}}.
\end{aligned}
\end{equation}
The integrals in (\ref{eq:NitscheFrom}) correspond to the standard consistency, adjoint consistency, and the Nitsche penalty terms, respectively. The third term is scaled by the Young's modulus, $E$, the element size, $h$, and the penalty factor $\gamma_N$. The latter provides additional control over the accuracy at which a boundary condition (BC) is enforced.

\subsubsection{Face-oriented ghost stabilization}

\label{sec:Ghost}
Face-oriented ghost penalization, as proposed by \cite{schott2014new} and \cite{burman2015cutfem}, is used in the vicinity of the interface. This stabilization technique prevents numerical instabilities due to vanishing zones of influence of certain DOFs occurring when the material interface moves too close to a node.
Independent of the intersection configuration, ill-conditioning is mitigated by applying the following virtual work-based formulation (\cite{geiss2018topology,geiss2019level}):
\begin{equation}\label{eq:GhostPenForm}
\begin{aligned}
	\mathbf{R}_\Gamma^G 
	= 
	h \gamma_G 
	\underbrace{\sum}_{\mathcal{F} \in \mathcal{F}_{cut}}
	\int_{\mathcal{F}}
	\llbracket \delta \nabla \mathbf{u} \mathbf{N_e} \rrbracket 
	\llbracket \boldsymbol{\sigma}(\mathbf{u}) \mathbf{N_e} \rrbracket 
	dA.
\end{aligned}
\end{equation}
Element faces in the vicinity of the material interface for which at least one of the two adjacent elements is intersected are included in $\mathcal{F}_{cut}$ (\cite{villanueva2017cutfem}).
Outward facing normals of these shared element faces are denoted by $\mathbf{N_e}$.
The influence of the ghost penalty term presented above is controlled by the penalty factor $\gamma_G$.
The virtual work-based formulation is adopted instead of the one proposed by \cite{burman2014fictitious} since it allows for different material properties in adjacent intersected elements.
In this work, the ghost penalty term is computed based on elementally constant material properties that are evaluated at element centroids (un-intersected elements) or sub-phase centroids (intersected elements). 

\subsubsection{Selective structural springs}
To suppress rigid body motion associated with disconnected material domains that may emerge and develop, selective structural springs as proposed in \cite{villanueva2017cutfem} are added to the governing equations. An additional stiffness is applied only to solid disconnected subdomains via the following residual component:
\begin{equation}\label{eq:SelecSprEq}
\begin{aligned}
	\mathbf{R}^S = \gamma_S ~ E / h^2 \int_{\Gamma_{D}}  \delta \mathbf{u} \cdot \mathbf{u} dV.
\end{aligned}
\end{equation}
The parameter $\gamma_S$ denotes the spring stiffness scaling and is non-zero only for the free-floating pieces of material. An auxiliary indicator field obtained by solving a diffusion problem is employed to identify such disconnected subdomains.
More details of the application of this approach to structural problems can be found in \cite{geiss2018topology,geiss2019combined}.

\subsubsection{Evaluation of stresses}
In the final example in Section \ref{sec:NumEx}, a gradient stabilized scalar stress field, $\tau(\mathbf{X})$, is post-processed via the XFEM informed smoothing procedure described in \cite{sharma2018stress}. The additional set of state variables is computed by solving the residual equation:
\begin{equation}\label{eq:stressProj}
\begin{aligned}
	\mathbf{R} =
	\mathbf{R}^{\tau}
	&= 
	\int_{\Omega_{I}} ~ \delta \tau ~ (\tau - \sigma_{VM}) dV \\
	&+
	h^2 \gamma_\tau
	\underbrace{\sum}_{\mathcal{F} \in \mathcal{F}_{cut}}
	\int_{\mathcal{F}}
	\llbracket \delta \nabla \tau \rrbracket 
	\llbracket \nabla \tau \rrbracket 
	dA.
\end{aligned}
\end{equation}
The field $\sigma_{VM}(\mathbf{X})$ and the parameter $\gamma_\tau$ represent the von Mises stress and the ghost penalty weight, respectively.
Overestimation of stresses is prevented by penalizing the jump in spatial stress gradients across elemental faces (second term in Eq. \ref{eq:stressProj}). 

\subsection{General optimization problem formulation} \label{subsec:TFCOptProbForm}

In this paper the following formulation of the optimization problem is considered: 
\begin{equation}\label{eq:TFCOptProbSetup}
\begin{aligned}
\underset{s}{\min}~z(\mathbf{s}, \mathbf{u (\mathbf{s})})  
& = 
 w_1 ~ \mathcal{F}(\mathbf{s}, \mathbf{u(\mathbf{s})})  
+ 
w_2  ~ P_{Per}(\mathbf{s})
\\
& 
+ 
w_3 ~ P_{Reg}(\mathbf{s})
+ 
w_4 ~ P_{\rho\phi}(\mathbf{s})
\\
s.t.: ~~~
&~g_i(\mathbf{s},\mathbf{u(\mathbf{s})}) \leq 0, i=1,...,N_g.
\end{aligned}
\end{equation}
The objective $z$ is minimized over the vector of admissible optimization variables, $\mathbf{s}$, defined in Sections \ref{subsec:SingleFieldCoupling} and \ref{subsec:TwoFieldCoupling} for the SFC and TFC, respectively; and the vector of state variables, $\mathbf{u}$, with $\mathbf{u} \in \rm I\!R^{N_{u}}$, and $N_u$ being the number of state variables.

The first component of the objective represents the quantity of interest to be minimized, $\mathcal{F}$ (e.g., strain energy, mass). 
The second term is the normalized perimeter control penalty, 
\begin{equation}\label{eq:PerPenForm}
\begin{aligned}
P_{Per}
= 
\frac{\displaystyle\int_{\Gamma_{I,II}} dA}{\displaystyle\int_{\Gamma_{D}} dA},
\end{aligned}
\end{equation}
which prevents the emergence of irregular geometric features. $\Gamma_{D}$ corresponds to the perimeter of the design domain $\Omega_{D}$.
The LS regularization penalty, $P_{Reg}$, is included in the objective to avoid spurious oscillations in the LSF (see Section \ref{subsec:LsReg}).
Note that $P_{Reg}$ promotes a positive LSF while the coupling penalty attempts to lower it to a negative value to nucleate a hole.
By reducing the weights $w_{\phi_2}$ and $w_{\phi\nabla2}$ in Eq. \ref{eq:LsRegWeights}, these competing effects are mitigated. This flexibility is exploited in the third example in Section \ref{sec:NumEx}.
The final component denotes the normalized two-field coupling penalty,
\begin{equation}\label{eq:couplingPenNormForm}
\begin{aligned}
P_{\rho\phi}
= 
\frac{\displaystyle\int_{\Omega_D^0} {p}_{\rho\phi}(\mathbf{X}) dV}{\displaystyle\int_{\Gamma_{D}}  dA},
\end{aligned}
\end{equation}
with ${p}_{\rho\phi}(\mathbf{X})$, as formulated in Eq. \ref{eq:TFCFomAllDomain}.
This coupling penalty can also be defined as a constraint. We choose to add it to the objective as it favors a looser coupling of the density and LS fields. In our experience, this approach is sufficient to nucleate holes and prevents over-constraining the design.

The weights $w_1$, $w_2$, $w_3$, and $w_4$ are chosen such that all the penalty contributions in the objective are significantly lower ($\approx 1-5\%$) than $\mathcal{F}$.
Since a strong perimeter control penalty might prevent the nucleation of small holes (\cite{wang2007hole}), its contribution is kept below $1\%$ by manipulating the $w_2$ weight. Both, a constant low weight and gradual increase of the perimeter contribution within the continuation scheme (i.e., $\mathcal{D}_{it} \leq \mathcal{D}_{c}$), are considered in the numerical examples.
In addition, for the remainder of the optimization process (i.e., $\mathcal{D}_{it}>\mathcal{D}_{c}$), the perimeter penalty weight is increased to promote a smooth geometry in the final design. 
For the SFC approach, $w_4$ is set to 0.0 since a coupling penalty is not required.
 
The design needs to satisfy a set of $N_g$ problem dependent inequality constraints, $[g_1,...,g_{N_g}]$ (e.g., target mass, maximum allowable stress, maximum eigenvalue).
The constraints are defined for each problem studied in Section \ref{sec:NumEx}.

\section{Numerical Examples}\label{sec:NumEx}

The proposed SFC and TFC approaches are studied in this section using single material, solid-void linear elastic problems in 2D and 3D.
Algorithmic parameter dependencies and performance are investigated with a 2D structural plane stress problem under uniform pressure. A beam problem in 3D is used to assess the influence of the SIMP penalization on the robustness of the coupling strategies. And finally, the overall behavior of these approaches is examined with a geometrically complex engineering problem considering stress constraints.

The optimization problems are solved using the Globally Convergent Method of Moving Asymptotes (GCMMA, \cite{svanberg2002class}) with no inner iterations. The adjoint method, as detailed in \cite{sharma2017shape}, is used for the sensitivity analysis. 
The relative change in objective between two consecutive design iterations is less or equal to $1\times10^{-3}$ at the end of all the problems shown in this section.
Relevant optimization parameters common to all problems presented in this work are summarized in Table \ref{tab:commonOptProbParams}.
The weights for the LS regularization penalty in Examples 1 and 2 are $w_{\nabla\phi_1}=w_{\phi_1}=w_{\phi_2}=w_{\nabla\phi_2}=1$, and for the third problem, $w_{\nabla\phi_1}=0$, $w_{\phi_1}=w_{\phi_2}=w_{\nabla\phi_2}=0.5$.

\begin{table}[h!] 
	\caption{\label{tab:commonOptProbParams}Parameters common to all numerical examples function of the element size $h$.}
	\center
	\renewcommand{\arraystretch}{1.2}
	\begin{tabular}{l|c}
		\hline
		Parameter                  & Value\\\hline		
LSF upper bound             		 &  $\phi_{up} = 2.5h$        \\
LSF lower bound             		        &  $\phi_{low} = -2.5h$       \\
Initial constant LSF       	    	        &  $\phi_{0} = 0.5\phi_{up}$ \\
LSF regularization control               &   $\gamma_{P_{Reg}} = 36.8$     \\
Filter radius  in 2D           	 &  $f_r = 1.6h $      \\
Filter radius in 3D           	 &  $f_r = 1.8h $      \\
Nitsche penalty factor &  $\gamma_N = 100.0 $      \\
Ghost penalty factor	 &  $\gamma_G = 0.001 $      \\
Spring stiffness factor           	 &  $\gamma_S = 1\text{x}10^{-6} $      \\
Stress ghost penalty factor (Ex. 3)          	 &  $\gamma_\tau = 0.01 $      \\
		\hline
	\end{tabular}
\end{table}

At the beginning of the optimization process, for both the SFC and TFC approaches, a uniform LSF of $\phi(\mathbf{X})=\phi_0$ $(\phi_0>0)$ is prescribed such that the entire design domain is filled with homogeneous porous material. 
In all examples, the continuation scheme described in Section \ref{subsec:DensTO} is used to update the shift in the densities unless specified otherwise; see Fig. \ref{fig:DensShiftScheme}. Furthermore, when using the TFC approach, the density threshold is decreased using the continuation scheme detailed in Section \ref{subsec:TwoFieldCoupling}; see Fig. \ref{fig:DensThresholdScheme}. 

The governing equations are discretized using the XFEM approach outlined in Section \ref{sec:XFEMFramework}. 
The parameters used for interpolation of material properties in Examples 1 and 2 are shown in Table \ref{tab:LsDensCombProp} in self-consistent units, and for the last example in Table \ref{tab:ex3BracketMatProps} in SI units.
The example problems consist of a one-way coupled set of governing equations; i.e., the diffusion problem describing the auxiliary indicator field, the stabilized linear elasticity equations (Eq. \ref{eq:resStrucGovEq}), the stabilized stress projection (Eq. \ref{eq:stressProj}), and the heat method equations.
Taking advantage of this simple coupling scheme, a Block Gauss-Seidel approach (\cite{elfving1980block}) is employed for both the physical and the adjoint sensitivity analyses.
The systems of equations of the first and second examples were solved using the Multifrontal Massively Parallel Solver (MUMPS, \cite{amestoy2006hybrid}).
The third example was solved iteratively via Trilinos algebraic multigrid solver (\cite{heroux2005overview}) for the selective springs heat problem, the structural response, and the stresses; and ILU preconditioner for the LS regularization heat method problem.

\begin{table}[h!]
	\caption{\label{tab:LsDensCombProp}Material properties interpolation parameters for numerical examples 1 and 2.}
	\center
	\renewcommand{\arraystretch}{1.2}
	\begin{tabular}{l|c}
		\hline
		Property                  & Value\\\hline
		Initial fictitious density  ($\Omega_{I}$)         & $\rho_0=0.4$              \\					
		SIMP exponent                                           &  $\beta_{\rho} $ = 2.0 \\			
		Young's Modulus ($\Omega_{I}$)  	           &  $E_S = 2.0$x$10^3$              \\
		Young's Modulus ($\Omega_{II}$)               &  $E_V = 1$x$10^{-8}$ \\
		Poisson Ratio ($\Omega_{I}$ and $\Omega_{II}$) &  $\nu_S = \nu_V$ = 0.4 \\			
		Material Density ($\Omega_{I}$)     	           &  $\theta_S = 1.0$ \\
		Material Density ($\Omega_{II}$)              	   &  $\theta_V = 0.0$ \\
		\hline
	\end{tabular}
\end{table}

\subsection{Example 1: structure subject to uniform pressure} \label{sec:SupportStruc2D}

\begin{figure}[h]
	\center
	\includegraphics[width=0.5\linewidth]{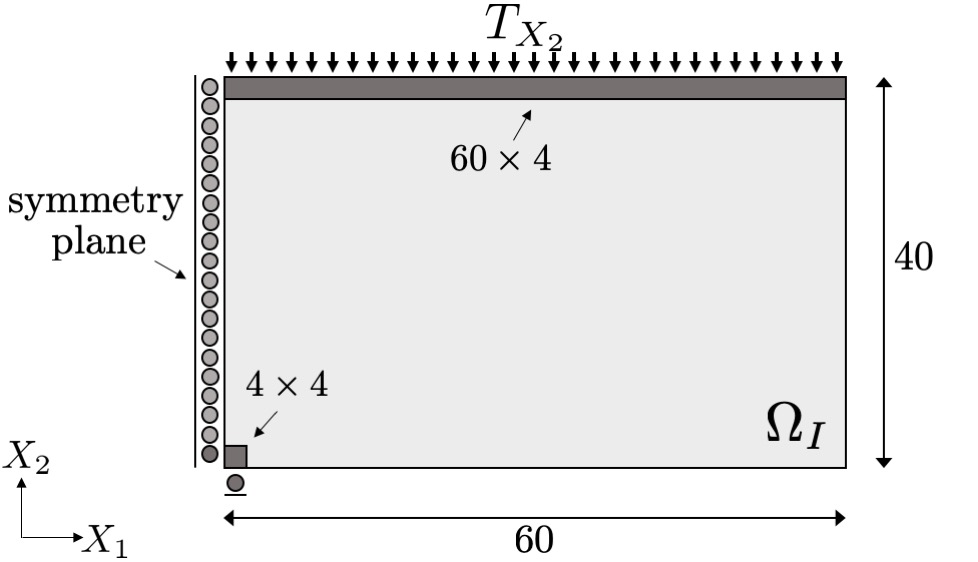}
	\caption{Problem setup of the 2D structural design problem. Due to symmetry, only half of the domain is simulated.}
	\label{fig:supportStruc2DProbSetup}
\end{figure}

A 2D structure subject to uniform pressure load is studied to demonstrate the hole nucleation capabilities of the proposed approaches. The problem setup with loads and BCs is shown in Fig. \ref{fig:supportStruc2DProbSetup}.
A vertical traction load of $T_{X_2} = -10.0$ is applied over the top boundary of the design domain. 
The bottom left corner is clamped by prescribing weakly zero displacement in the $X_1$ and $X_2$ directions.
Design symmetry about the $X_2$ axis is assumed in the middle of the design domain (i.e., zero displacement in the $X_1$ direction along symmetry plane). 
Solid material is prescribed and excluded from the design domain in the vicinity of the Dirichlet BCs and where the traction load is applied, as highlighted in dark grey color in Fig. \ref{fig:supportStruc2DProbSetup}.
The domain of size $60\times40$ is discretized by a structured mesh with uniform element size $h=0.5$.

The compliance minimization with mass constraint optimization problem reads:
\begin{equation}\label{eq:Ex1OptProbSetup}
\begin{aligned}
\underset{s}{\min}~z(\mathbf{s}) & = 
w_1~ \Psi (\mathbf{s},\mathbf{u(\mathbf{s})}) /\Psi_0 
+ 
w_2 ~ P_{Per}(\mathbf{s})
\\
& 
+ 
w_3 ~ P_{Reg}(\mathbf{s})
+ 
w_4 ~ P_{\rho\phi}(\mathbf{s})
\\
s.t.: g_1 & = \frac{\mathcal{M(\mathbf{s})}}{\Omega_{I}+\Omega_{II}} - \gamma_m \leq 0.
\end{aligned}
\end{equation}
The strain energy of the initial design is denoted by $\Psi_0$, and the objective weights for the TFC setup are $w_i$ =[ 0.93, 0.01, 0.05, 0.01]. In the SFC setup, $w_4=0$. 
A mass constraint of $\gamma_m = 0.40$ is enforced.

The default continuation parameters for updating the density shift (in both approaches) and the density threshold (in the TFC approach) are summarized in Table \ref{tab:Ex1DefContParams}, and used for all example configurations unless specified otherwise.
Both $\rho_{sh}$ and $\rho_{th}$ are updated every $50$ optimization iterations in a total of eight continuation steps. 
No density shift is applied at the beginning of the optimization process. All simulations start with a uniform density of $\rho_{0}=0.40$, except the portions excluded from the design domain, which are prescribed a fixed density of 1.0.
Note that the SIMP exponent $\beta_\rho$ is set to 2.0 to mitigate mesh dependency, avoid premature convergence, and allow for nucleating many holes. An analysis of the influence of this parameter is provided in Example 2.

\begin{table}[h!]
	\caption{\label{tab:Ex1DefContParams}Default continuation parameters for the 2D structural problem.}
	\center
	\renewcommand{\arraystretch}{1.2}
	\begin{tabular}{l|c}
		\hline
		Parameter                  & Value\\\hline
		Continuation step size&  $\mathcal{D}_{st}=50$              \\
		Number of design iterations in continuation &  $\mathcal{D}_{c}=400$              \\		
		Maximum number of design iterations &  $\mathcal{D}_{max}=500$              \\			
		Initial density shift &  $\rho^0_{sh}=0.0$              \\			
		Initial density threshold &  $\rho^0_{th}=0.7 \rho_0$              \\
		Continuation density shift exponent &  $\eta_{\rho_{sh}}=2.0$              \\
		Continuation density threshold exponent &  $\eta_{\rho_{th}}=2.0$       \\		
		\hline
	\end{tabular}
\end{table}

\subsubsection{Influence of density shift} \label{subsubsec:SupportStruc2DIntermediateDens}

\begin{figure}[h]
	\center
	\includegraphics[width=0.5\linewidth]{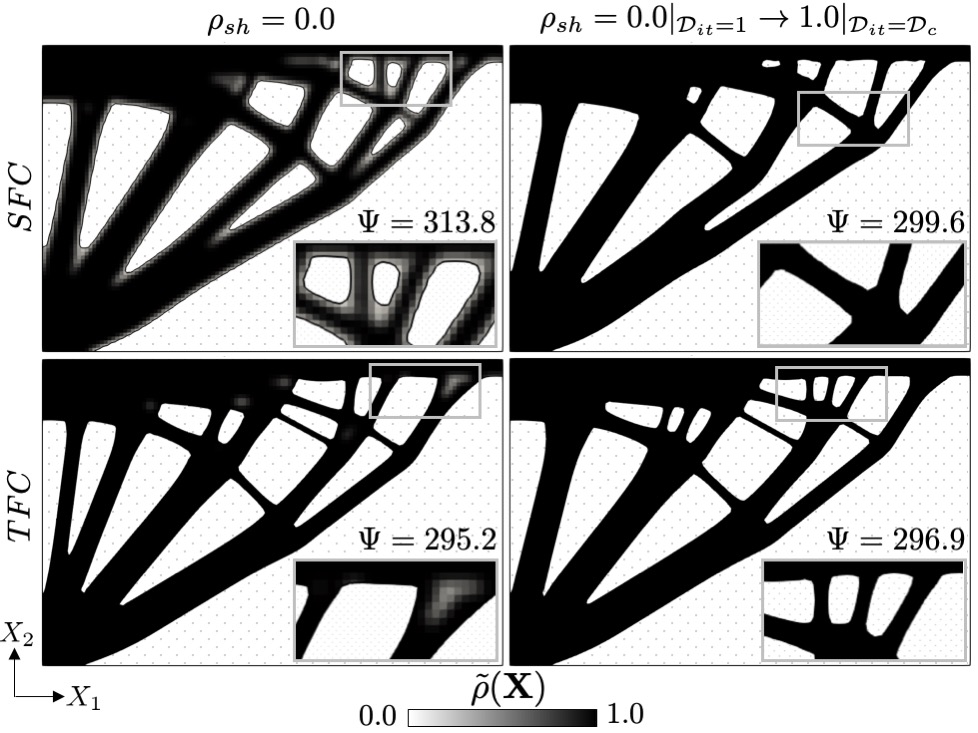}
	\caption{Final designs of the 2D structural problem for both SFC and TFC approaches with and without the augmented SIMP material interpolation.}
	\label{fig:supportStructShiftAndNoShiftResults}
\end{figure}

In this first study, the effect of the density shift described in Section \ref{subsec:DensTO} is investigated. 
The final designs for both the single and two-field coupling approaches with and without density shift are shown in Fig. \ref{fig:supportStructShiftAndNoShiftResults}. 
Only the material phase, $\Omega_{I}$, is visible; and it is colored by the shifted filtered densities, $\tilde\rho(\mathbf{X})$. The void regions are identified by dotted patterns with white background.
The SFC and TFC results are shown in rows one and two, respectively, without and with a density shift (left and right columns). 
The insets highlight regions of the final designs where significant changes in the density distributions are observed.

On the left side of Fig. \ref{fig:supportStructShiftAndNoShiftResults} it can be seen that both methods are unable to fully remove intermediate densities in absence of a density shift. In the SFC approach, since both the LSFs and fictitious densities are function of the same abstract design variables (Eqs. \ref{eq:SFC_AdvToLsVars} and \ref{eq:SFC_AdvToDensVars}), intermediate densities exist in the vicinity of the material interface by construction. In the TFC, despite the density field can evolve independently of the LSF, grey areas are also observed.

The results in the right column of Fig. \ref{fig:supportStructShiftAndNoShiftResults} show that the gradual shift scheme adopted guarantees a uniform density distribution of 1 in the material subdomain.
Moreover, even though the topologies being noticeably different, the performance of the final design is not compromised by shifting the densities. On the contrary, a considerable improvement is observed in the SFC results.
Note that as an alternative to improve the convergence to [0-1] designs, the SIMP penalization could be increased and combined with a projection scheme (e.g., see \cite{wang2011projection}). However, those approaches do not guarantee that the solid phase consists only of bulk material in the optimized design. See Example 2 for a more detailed discussion on this matter.

~\\
\textit{Initial density shift}

To better understand the influence of an initial density shift in the hole nucleation process, a comparison of the design evolutions with $\rho^0_{sh}>0.0$ is presented. All parameters but the initial density shift are the same as before. Given that the numerical examples in this study show similar findings for both the SFC and TFC approaches, we only show results for the SFC approach.
\begin{figure*}[h]
	\center
	\includegraphics[width=1.0\linewidth]{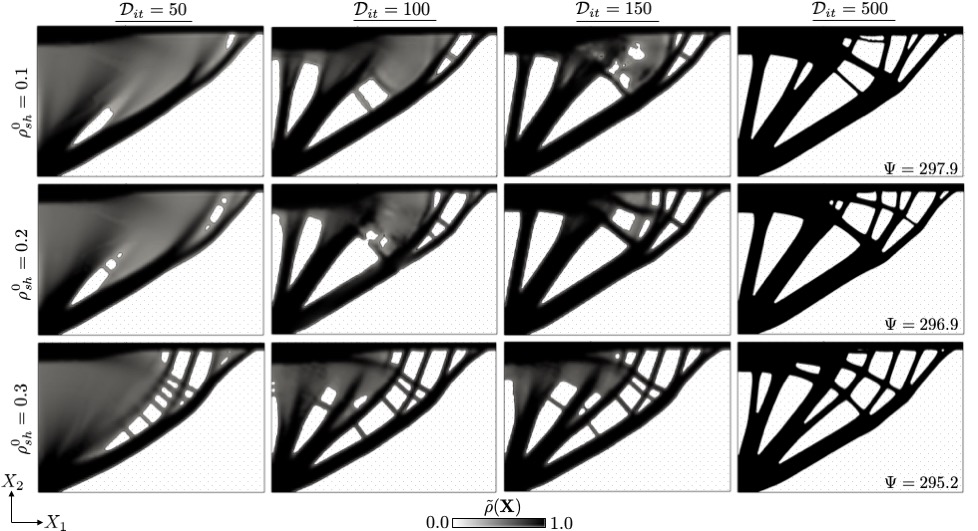}
	\caption{Snapshots of the evolution of the 2D structural design using the SFC approach from left to right are shown with different initial density shift per row.}
	\label{fig:influenceOfShiftSFC}
\end{figure*}

Snapshots of the optimization process at $\mathcal{D}_{it}$ = [50, 100, 150, 500] for $\rho_{sh}^0$ = [0.1, 0.2, 0.3] are provided in Fig. \ref{fig:influenceOfShiftSFC}. A fast and clean hole nucleation process is observed in all cases. 
Note that densities in the void phase do not contribute to neither the stiffness nor the mass of the design. Therefore, sensitivities associated with these variables only exist in the material domain ($\Omega_{I}$).
Material is removed from the design domain by either nucleating a hole due to low densities or shrinking/expanding the material interface front.
A larger initial density shift promotes the generation of more holes at an early stage, as seen in the first column of Fig. \ref{fig:influenceOfShiftSFC}. Nucleating less holes initially (guided by the density sensitivities) is, however, balanced by the moving material interface (result of shape sensitivities of the LS problem). Eventually, similar topologies at a later stage of the optimization process are observed. Although most holes are nucleated within the first third of the optimization process, the hole nucleation process is active, and holes are nucleated until the density shift enforces the maximum density at $\mathcal{D}_{it}=\mathcal{D}_{c}$. 

The performance of the final design is not compromised by introducing a density shift. As is reported on the last column of Fig. \ref{fig:influenceOfShiftSFC}, the variation of the final strain energy is less than 1\%. Considering the results of the previous subsection where $\rho^0_{sh}=0.0$ (see Fig. \ref{fig:supportStructShiftAndNoShiftResults}), there is a small improvement on performance with an increasing shift.

\subsubsection{Influence of continuation strategy} \label{subsec:SupportStruc2DDensShift}

\begin{figure}[h!]
	\center
	\includegraphics[width=0.5\linewidth]{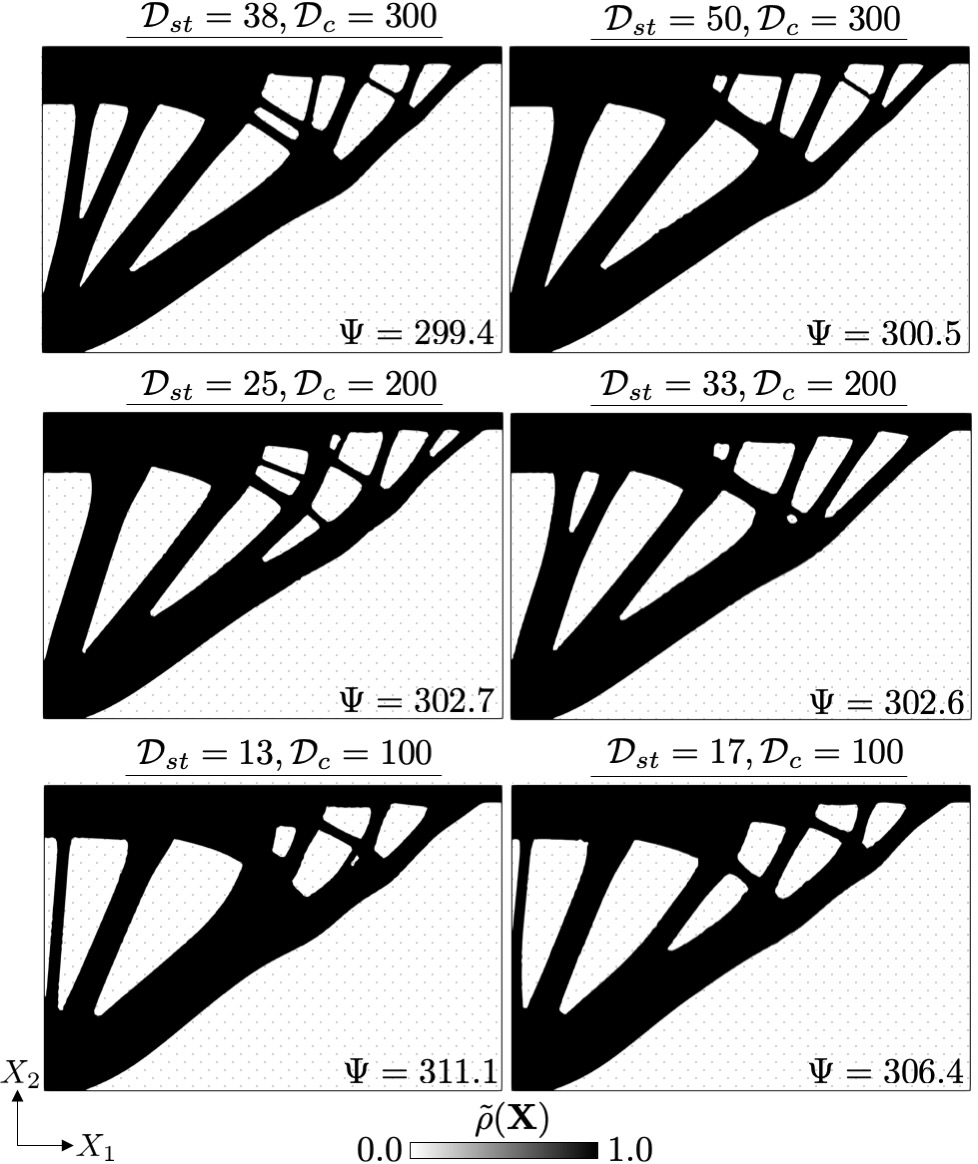}
	\caption{Final designs of 2D structural problem for multiple continuation step sizes and number of design iterations with a continuation strategy.}
	\label{fig:influenceOfContStepSizeAndSpam}
\end{figure}

\begin{figure}[h!]
	\center
	\includegraphics[width=0.5\linewidth]{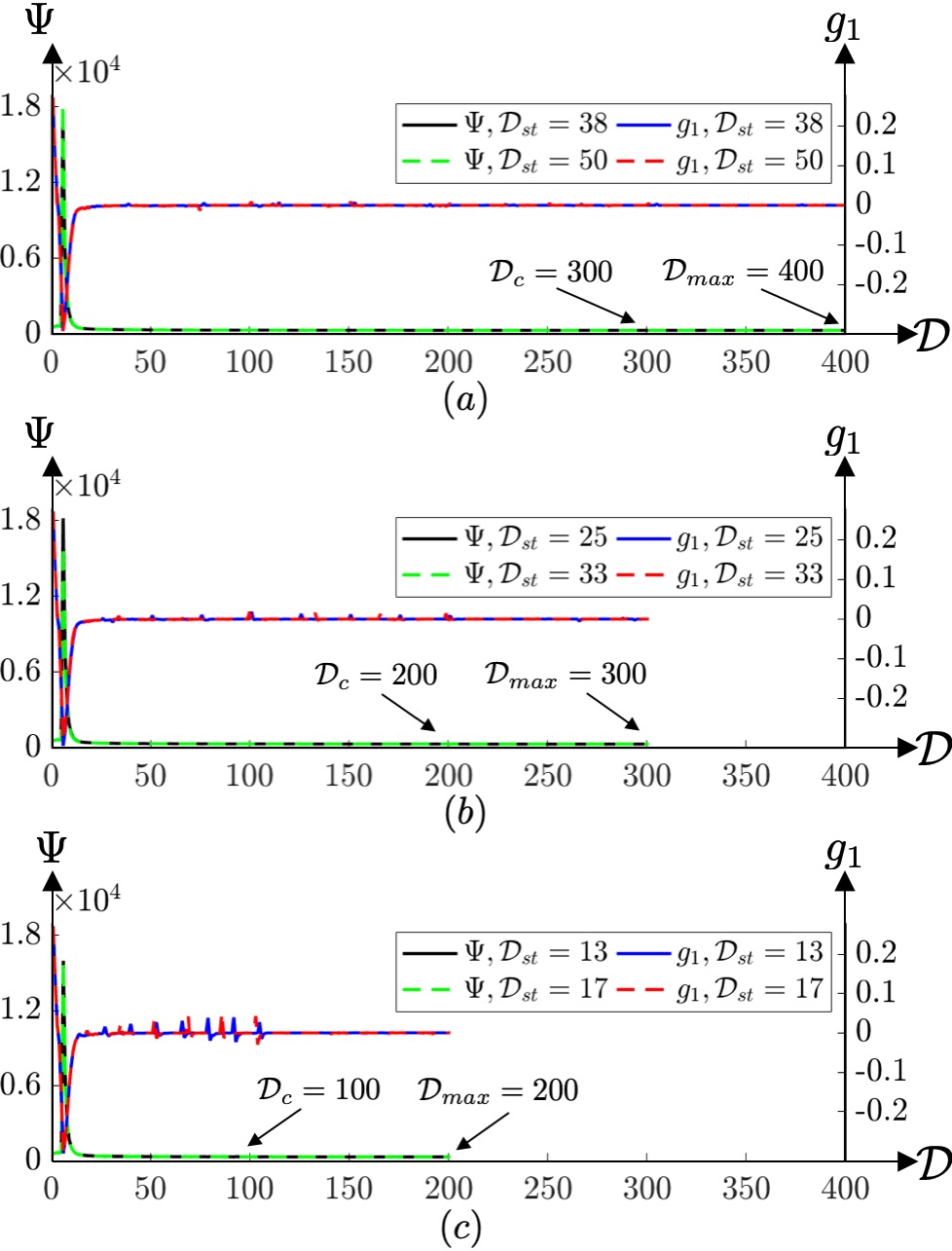}
	\caption{Evolution of objective and constraint of the 2D structural problem for two continuation step sizes and multiple total number of design iterations.}
	\label{fig:SFC_ObjAndConstEvol_contStepAndRange}
\end{figure}

The final designs of different continuation settings for the density shift using the SFC approach are shown in Fig. \ref{fig:influenceOfContStepSizeAndSpam}. 
A total of six (left) and eight (right) continuation steps are used in optimization problems with (a) 400, (b) 300, and (c) 200 design iterations. 
The continuation scheme parameters for each case are provided on top of the final designs.
These converged designs demonstrate that the proposed approaches provide the flexibility of choosing different continuation step sizes, $\mathcal{D}_{st}$, and total number of continuation iterations, $\mathcal{D}_{c}$, to accelerate the design process. Similar performances are achieved for all configurations.
However, reducing the continuation step size decreases the number of holes nucleated. An excessively small continuation step size may cause the optimization process to prematurely converge to a suboptimal design. 
Fig.\ref{fig:SFC_ObjAndConstEvol_contStepAndRange} shows the evolution of the objective and constraints of the designs in Fig. \ref{fig:influenceOfContStepSizeAndSpam}, separated by total number of design iterations.
Here, $\mathcal{D}_{max}$ represent the maximum number of design iterations allowed.
It can be seen that a larger number of continuation steps smoothes the transition of the density shift at the cost of more frequent modifications to the optimization problem. 

Numerical experiments with the problems studied in this paper suggest that using five or more continuation steps of $\mathcal{D}_{st}\approx[15,25]$ provides an adequate seeding capabilities and circumvents large jumps in the densities updates.

\subsubsection{TFC optimization problem characteristics} \label{subsubsec:SupportStruc2DObjAndSensTFC}

\begin{figure}[h!]
	\center
	\includegraphics[width=0.5\linewidth]{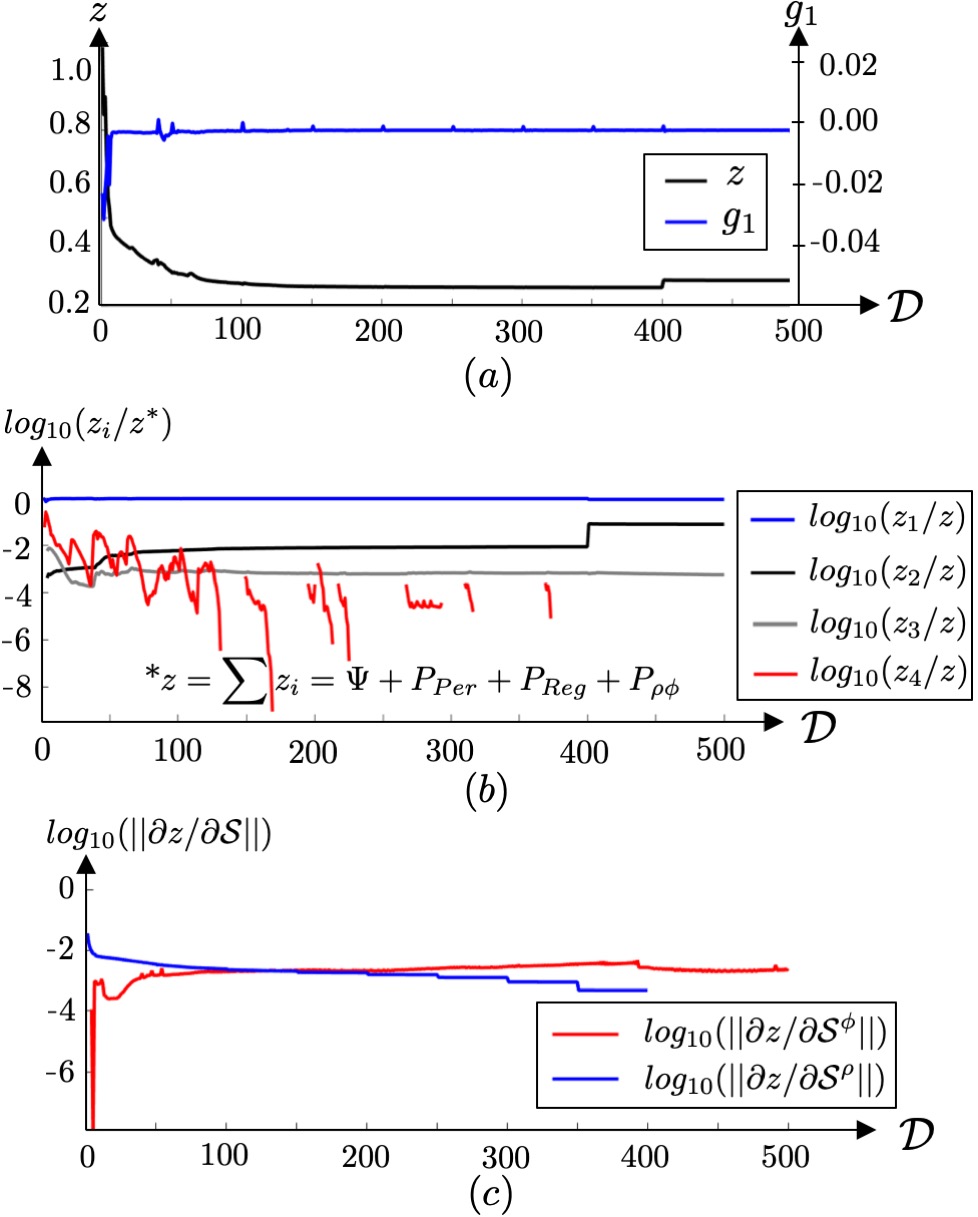}
	\caption{Evolution of (a) objective and mass constraint, (b) logarithm of contributions of the components of the objectives; and (b) implicit sensitivities contributions (logarithms of norms) of the density and LS design variables for the 2D structural problem using the TFC approach.}
	\label{fig:TFCObjAndConst_CompContr_SensContr}
\end{figure}

Fig. \ref{fig:TFCObjAndConst_CompContr_SensContr} shows (a) the evolution of the objective and mass constraint, (b) the objective components contributions, and (c) the LS and density sensitivity contributions to the optimization problem for the TFC results with density shift of the 2D structural problem in Fig. \ref{fig:supportStructShiftAndNoShiftResults}.
Small oscillations are observed early on in the design process due to holes being nucleated; see Fig. \ref{fig:TFCObjAndConst_CompContr_SensContr}(a).
Since this occurs in low density regions, the nucleation process has a reduced effect on the objective and thus, the observed oscillations are small and quickly vanish.
A smooth design evolution is achieved in a later stage. This is in contrast to using alternatives like topological derivatives for which, in our experience, the size and number of holes nucleated may produce large jumps in the optimization process.
A slight increase in the objective at $\mathcal{D}_{c}=400$ originates from increasing the weight of the perimeter control penalty. However, its impact on the strain energy convergence is negligible.
Small jumps in the mass constraint, $g_1$, at every continuation step update are consequence of sudden increments in mass due to the density shift.

~\\
\textit{Objective components contributions}

The logarithm of the contributions of each objective component (i.e., $\Psi$, $P_{Per}$, $P_{Reg}$, and $P_{\rho\phi}$; see Eq. \ref{eq:Ex1OptProbSetup}), are provided in Fig. \ref{fig:TFCObjAndConst_CompContr_SensContr}(b).
The optimization problem is at all times completely dominated by the strain energy.
The coupling, LS regularization and perimeter control penalties are on average two orders of magnitude lower. 
As a result, the evolution of the objective is predominantly smooth despite considerable oscillations and periods of inactivity of $P_{\rho\phi}$; see Fig. \ref{fig:TFCObjAndConst_CompContr_SensContr}(b).
The intermittent non-zero TFC penalty contribution is a direct consequence of its definition in Eq. \ref{eq:TFCFomAllDomain}. Intermediate densities larger than $\rho_{th}$ may shift up or down, when a hole is nucleated or due to the moving interface front. 
Sudden increments in the coupling penalty represent an intermediate density dropping below $\phi_{th}$. 

~\\
\textit{LS and Density sensitivities}

The interplay of the LS and density variables in the TFC approach is examined with Fig. \ref{fig:TFCObjAndConst_CompContr_SensContr}(c). The logarithm of the norms of the sensitivities with respect to the LS and density variables show that initially the density field drives the optimization process. Then, gradually with the introduction of holes, the shape sensitivities take over. The influence of the density sensitivities vanishes once the shifted densities have reached the maximum density at $\mathcal{D}_{c} = 400$.
A pure LS problem with shape sensitivities only is recovered for $\mathcal{D}_{it} > \mathcal{D}_{c}$.

\subsubsection{TFC thresholds}

As explained in Section \ref{subsec:TwoFieldCoupling}, the TFC coupling penalty is bounded by the LS ($\phi_{th}$) and the density ($\rho_{th}$) thresholds.
To avoid robustness issues in the hole nucleation process, $\phi_{th}$ is set to a number below the zero isocontour ($\phi_{th}= 0.25\phi_{low}$ in this example). As long as this condition is satisfied, no significant changes are observed in the hole nucleation process.
Similarly, the density threshold is constrained by the uniform density field $\rho(\mathbf{X})=\rho_0$ used at the beginning of the optimization process. 
However, unlike the LS threshold, different settings of the density threshold can have a considerable impact in the initial stages of the hole nucleation process, as is demonstrated below.

\begin{figure*}[h!]
	\center
	\includegraphics[width=1.0\linewidth]{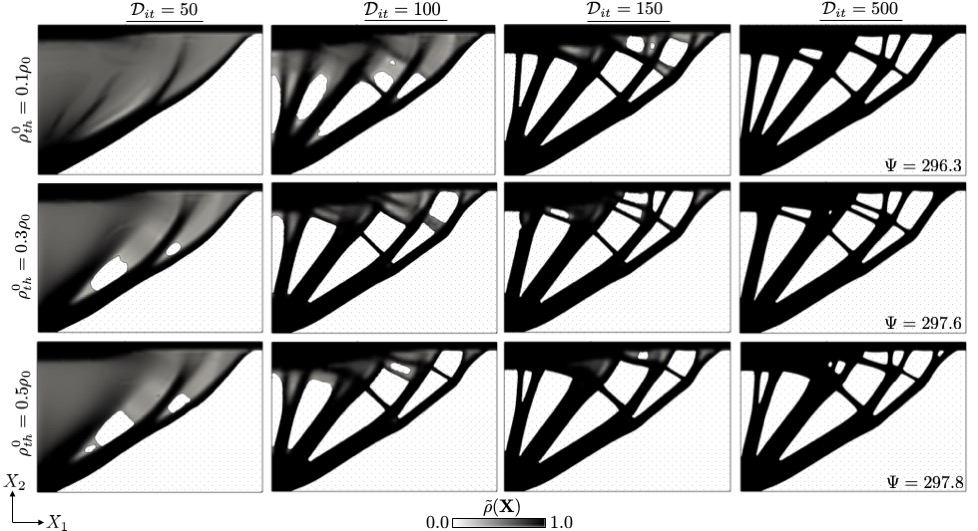}
	\caption{Snapshots of the evolution of the 2D structural problem using the TFC approach from left to right are shown using different initial density threshold per row.}
	\label{fig:influenceOfDensThTFC}
\end{figure*}

Snapshots of the evolution of the 2D structural design problem are presented in Fig. \ref{fig:influenceOfDensThTFC} to assess the influence of setting the initial density threshold closer to the initial constant density field.
The designs at $\mathcal{D}_{it}$=[50, 100, 150, 500] for $\rho_{th}^0$=[$0.1\rho_0$, $0.3\rho_0$, $0.5\rho_0$] separated by rows, show significant variations in the hole nucleation process.
As expected, setting $\rho^0_{th}$ closer to $\rho_0$ promotes the nucleation of holes earlier given that the penalty is active earlier in the design process.
The final designs for different $\rho_{th}^0$ values exhibit slightly different topologies but similar performances.
Also, as seen for the SFC approach previously, there is an interplay between the densities and the LS interface front until the density shift removes all intermediate densities. As long as intermediate densities exist, a hole can be nucleated if it is advantageous.

Overall, similar performances are achieved using a wide range of choices for parameters in both the single- and two-field strategies, demonstrating the robustness of these approaches with respect to algorithmic parameters.

\subsection{Example 2: beam} \label{sec:MbbBeam3D}

\begin{figure}[h!]
	\center
	\includegraphics[width=0.5\linewidth]{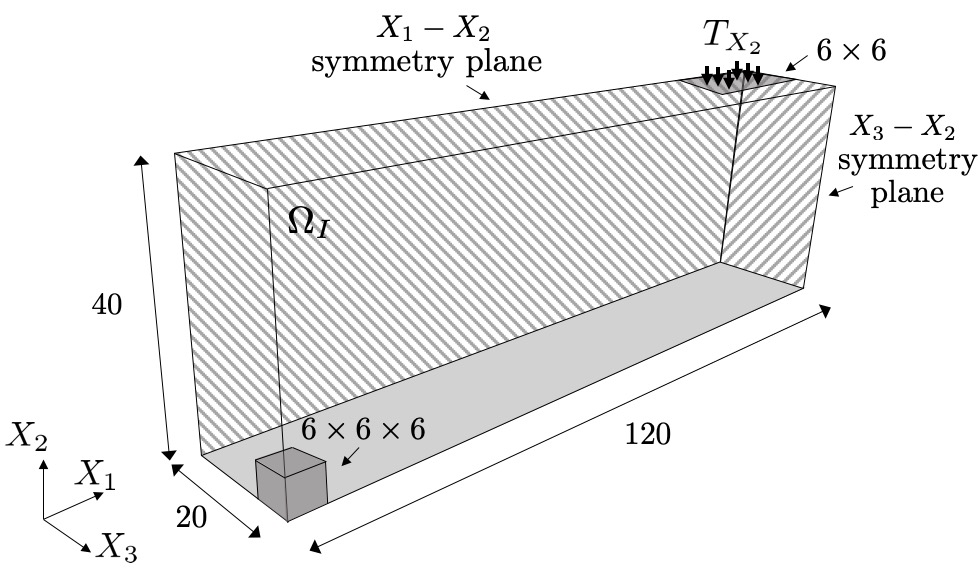}
	\caption{Problem setup for the classical Beam problem in 3D. Due to symmetry, only one quarter of the domain is simulated.}
	\label{fig:3dMbbBeamProbSetup}
\end{figure}

A 3D beam problem is used to investigate the effect of the SIMP penalizations in the density field, as well as the behavior of both approaches compared to LS-based with initial hole seeding and density-based TO.
The design domain of size $240\times40\times40$ is simply supported weakly at all four corners at the bottom face.
A traction load $T_{X_2} = -10.0$ is applied over the center of the top face.
Due to symmetry, only one quarter of the full domain is analyzed and optimized. 
The problem setup is shown in Fig. \ref{fig:3dMbbBeamProbSetup}.

The optimization problem formulation remains the same as in the previous example (see Eq. \ref{eq:Ex1OptProbSetup}). The initial objective weights are $w_i$ = [0.92, 0.001, 0.01, 0.05]. The perimeter weight, $w_2$, is updated at every continuation step from 0.001 to 0.01, using Eq. \ref{eq:expFuncContInc} with an exponent of 3.0. A mass constraint of $\gamma_m = 0.20$ is enforced.

The density field is shifted using the augmented SIMP continuation scheme discussed in Section \ref{subsec:DensTO}; and the density threshold for the TFC approach is gradually reduced to zero using Eq. \ref{eq:expFuncContDec}.
The continuation parameters are listed in Table \ref{tab:Ex2DefContParams}. The design domain is initialized with a uniform density of $\rho_0=0.2$ everywhere but along the prescribed BCs, where the density is set to 1.0 for the duration of the optimization process.

\begin{table}[hb!]
	\caption{\label{tab:Ex2DefContParams}Continuation parameters for the beam problem.}
	\center
	\renewcommand{\arraystretch}{1.2}
	\begin{tabular}{l|c}
		\hline
		Parameter                  & Value\\\hline
		Continuation step size&  $\mathcal{D}_{st}=20$              \\
		Number of design iterations in continuation &  $\mathcal{D}_{c}=120$              \\		
		Maximum number of design iterations &  $\mathcal{D}_{max}=150$              \\			
		Continuation density threshold exponent &  $\eta_{\rho_{th}}=2.0$       \\
		Continuation density shift exponent &  $\eta_{\rho_{sh}}=2.0$              \\
		Initial density threshold &  $\rho^0_{th}=0.75 \rho_0$              \\
		Initial density shift &  $\rho^0_{sh}=0.2$              \\			
		\hline
	\end{tabular}
\end{table}

\subsubsection{SIMP density penalization} 

Large exponents (i.e., $\beta_{\rho} \geq 3.0$) reduce intermediate densities at the cost of more nonlinear optimization problems and poorer conditioned finite element problems. 
Contrary to these methods, to converge to a [0-1] density field we employ the density shift scheme presented in Section \ref{subsec:DensTO}. 
To examine the reliance of the proposed approaches on the penalization of the densities for the material description (see Eq. \ref{eq:varyingMatPropForm}), different exponents in the power law are considered.  

\begin{figure}[h!]
	\center
	\includegraphics[width=0.5\linewidth]{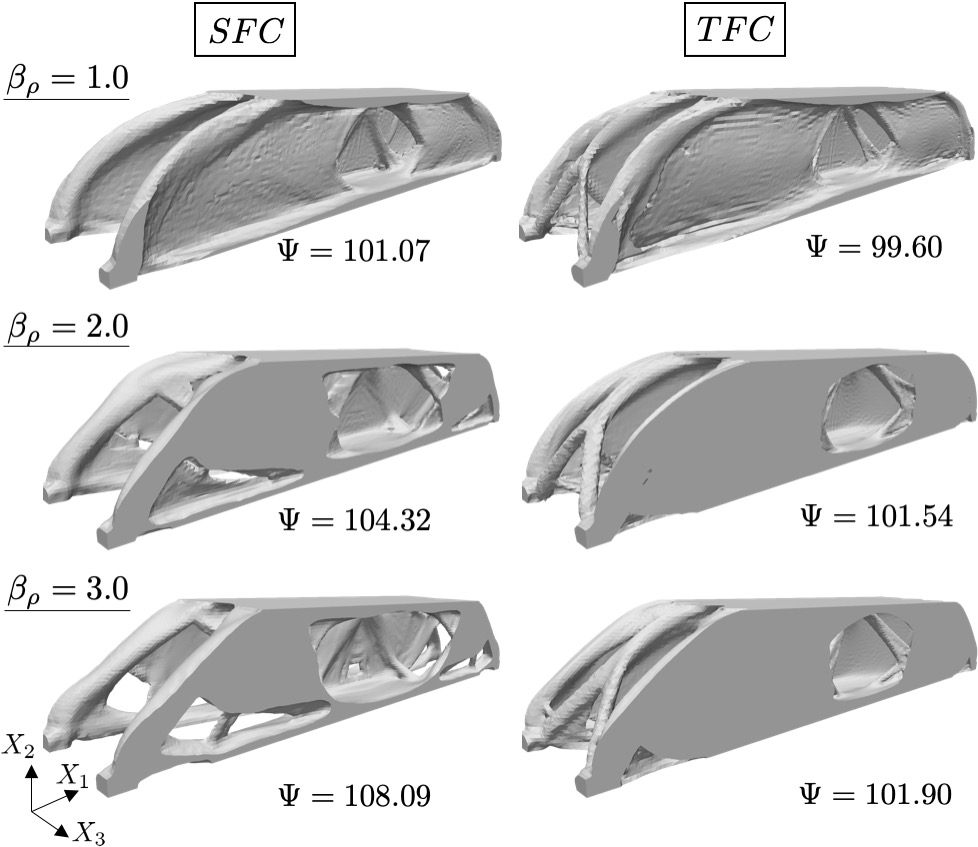}
	\caption{Final designs for both the SFC (left), and TFC (right) approaches using SIMP exponents of $\beta_{\rho}=[1.0, 2.0, 3.0]$.}
	\label{fig:mbbBeamFinalCombinedDesigns}
\end{figure}
%

\begin{figure*}[h!]
	\center
	\includegraphics[width=1.0\linewidth]{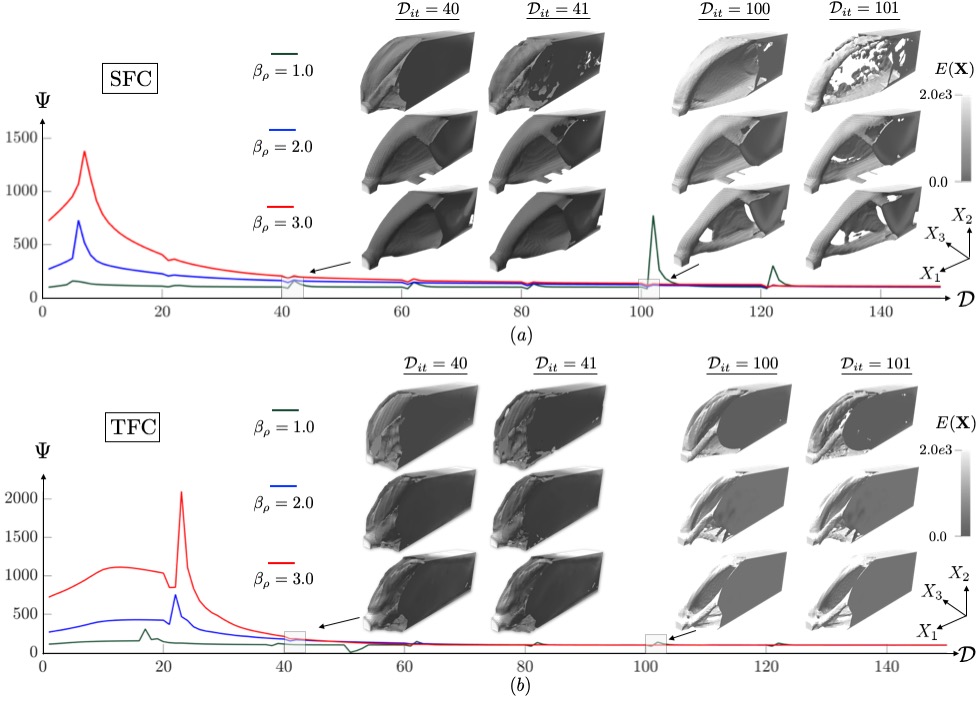}
	\caption{Evolution of the strain energy in the 3D beam problem using the (a) SFC  and (b) TFC approaches. Snapshots of designs before and after updating continuation parameters are colored by the Young's modulus.}
	\label{fig:mbbBeamSimpExpEvolSFCAndTFC}
\end{figure*}

The final designs of the beam problem using the single and two-field strategies with penalizations of $\beta_{\rho}$=[1.0, 2.0, 3.0] are shown in Fig. \ref{fig:mbbBeamFinalCombinedDesigns}. Their strain energy is also reported to evaluate structural performance.
All designs converge to a [0-1] density distribution as a consequence of the density shift. 
Note that, although not shown, intermediate densities were also observed without the shifting scheme, in accordance with the findings in Example 1 for the 2D problem.
For all $\beta_{\rho}$ values, similar truss-like structures with partial shear-webs are obtained. When using the TFC approach, the performances for the designs does not differ significantly with changing $\beta_{\rho}$. For the SFC approach, however, the strain energy of the optimized design increases noticeably with increasing $\beta_{\rho}$ value. Due to the tight coupling for density and LS fields, and the penalizations of intermediate densities, a more truss-like structure with inferior performance is developed as the density penalization is increased.

Fig. \ref{fig:mbbBeamSimpExpEvolSFCAndTFC} contains the objective evolution using both approaches together with snapshots of the quarter design domain colored by the Young's modulus to visualize the effect of the penalization. Designs before and after updating the shifted density field with a continuation scheme are presented for all $\beta_{\rho}$. A considerable initial increase in the strain energy, $\Psi$, is observed in all cases at early stages of the optimization process. The magnitude of the peaks depends on the strength of the density penalization since more material with intermediate densities is removed as $\beta_{\rho}$ increases.
A small penalization slows the hole nucleation process.
In the TFC approach, a low initial density threshold $\rho_{th}^0$ further delays hole nucleation; see also Section \ref{subsubsec:SupportStruc2DObjAndSensTFC}. 

We found that increasing the density shift parameter, $\rho_{sh}$, in the continuation strategy can compromise robustness if a small $\beta_{\rho}$ is employed, as seen in Fig. \ref{fig:mbbBeamSimpExpEvolSFCAndTFC}(a) for the SFC approach.
The evolution of the strain energy is mostly smooth with small influences of the $\beta_{\rho}$ parameter for $\beta_{\rho}$ = [2.0, 3.0]. 
However, due to a more diffuse density field with intermediate densities, applying the density shift can result in abrupt design changes, as seen for $\beta_{\rho}=1.0$. 
Since the shift increments become larger with every continuation step if $\eta_{\rho_{sh}}>1$ (see Fig. \ref{fig:DensShiftScheme}), this effect is more pronounced later in the optimization process. This can be especially detrimental in designs that favor thin features and for the SFC approach where intermediate densities are found in the vicinity of the interface by construction. The snapshots of the designs at $\mathcal{D}_{it}=100$ and $\mathcal{D}_{it}=101$ for $\beta_{\rho}=1.0$ in Fig. \ref{fig:mbbBeamSimpExpEvolSFCAndTFC}(a) illustrate this scenario. 
In contrast, a less noticeable impact of updating the density shift on the design evolution is observed in the TFC approach, as seen in Fig.\ref{fig:mbbBeamSimpExpEvolSFCAndTFC}(b). In this case the densities transition to 1.0 faster, especially near the interface. 

Overall, the stability and robustness of both coupling methods rely to a certain extent on penalizing the intermediate densities. For the problems studied in this paper, values for SIMP penalization in the range of $\beta_{\rho}$=[2.0-3.0] allow for the formation of fine geometric features and a smoother optimization process while avoiding ill-conditioning. 
For other problems where a SIMP penalization is insufficient, alternative techniques for removing intermediate densities might need to be considered.

\subsubsection{Comparing TO approaches} 

\begin{figure*}[h]
	\center
	\includegraphics[width=1.0\linewidth]{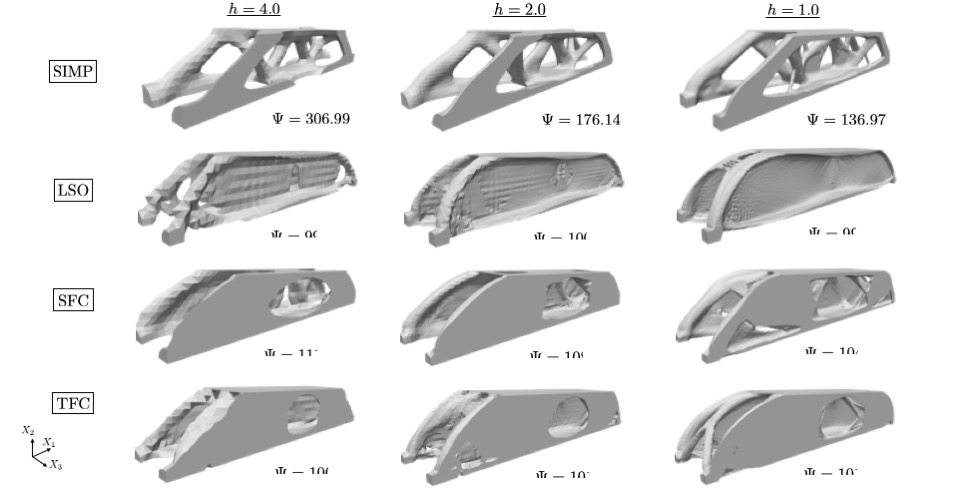}
	\caption{Final designs using meshes with element sizes $h$=[4.0, 2.0,1.0] of classical SIMP, LS-XFEM only (LSO) with initial hole seeding, and using the SFC and TFC approaches.}
	\label{fig:mbbBeamLsAndDensOnlyAndCombAppMeshRef}
\end{figure*}

The mesh dependency and performance of the SFC and TFC in comparison to LS TO with initial hole seeding, LSO (LS only), and standard SIMP are investigated. 
Final designs obtained by all four TO approaches are shown Fig. \ref{fig:mbbBeamLsAndDensOnlyAndCombAppMeshRef} for three levels of refinement, i.e., element lengths of $h=[4.0, 2.0, 1.0]$. In the LSO setups, 12 holes were initially seeded and a total of 300 design iterations without continuation were employed. 
The final SIMP designs are post-processed by applying an isovolume filter on the density field such that a total volume of 1.2$\times10^4$ is preserved in all designs.
Filtered densities are penalized by an exponent of 2.0 in the standard SIMP, SFC and TFC results.
The same filter radius is used in all cases; see Table \ref{tab:commonOptProbParams}. 

%

\begin{figure*}[h]
	\center
	\includegraphics[width=1.0\linewidth]{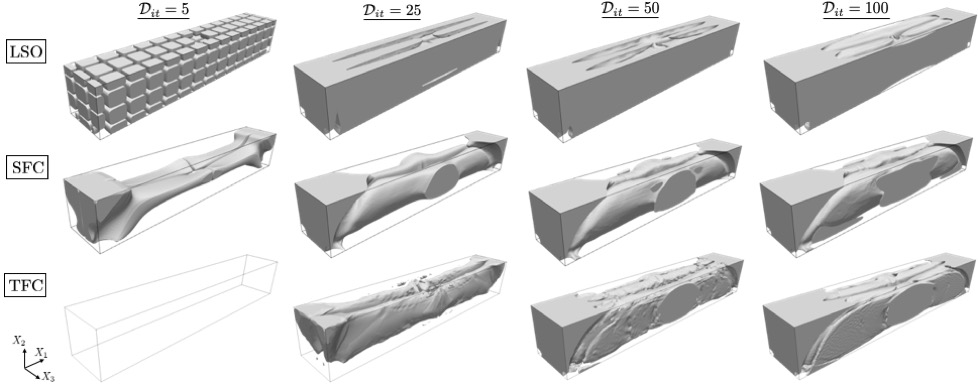}
	\caption{Snapshots of void domain for LSO with initial hole seeding, and using the SFC and TFC approaches for $h=1.0$.}
	\label{fig:mbbBeamLsApproachesVoidEvolComparison}
\end{figure*}

The SIMP results transition from a truss-like design to a shear-web structure as the mesh is refined (see first row of Fig. \ref{fig:mbbBeamLsAndDensOnlyAndCombAppMeshRef}). Suboptimal performances are attained in these designs due to intermediate densities not completely removed with a constant SIMP exponent of 2.0. The coarser the mesh, the more noticeable is this effect. 
A SIMP approach that projects the density field (e.g., see \cite{wang2011projection}) could have been used instead to mitigate the intermediate densities issue. However, to keep consistency in the comparison and assess the resemblance of the density problem and the coupling approaches in terms of their topologies, the same density penalization scheme is used in all cases.

The LSO final designs are continuous walls connecting the two flanges at the top and bottom of the beam with a thickness of about the element edge length.
Wrinkles are observed in all three levels of refinement.
However, they are reduced as the mesh is refined.
These spatial oscillations in the topology are typically mitigated by using either a larger filter radius or higher order interpolation functions. 
Nevertheless, in order to keep consistency in the comparison with the proposed approaches, the same filter radius and linearly interpolated fields are used in all runs. 

As reported in the last two rows of Fig. \ref{fig:mbbBeamLsAndDensOnlyAndCombAppMeshRef}, the SFC and TFC approaches experience smaller variations in final topologies and performance with mesh refinement. 
In both cases, the density field promotes final designs with shear-webs that are typically seen for LS-XFEM results with initial hole seeding.
Despite using a density field, they retain the desired LS-XFEM feature of having an increased geometrical resolution in design with coarser meshes. Furthermore, the TFC approach seems to provide more consistent results with mesh refinement.
Note also that the SFC and TFC strategies have an inherent tendency to form bulkier features due to the SIMP penalization.
In settings where the initial hole pattern avoids triggering suboptimal final designs, like the one shown in this example, pure LS results are marginally better than using the SFC and TFC strategies. However, this is might not the case in general.

Fig. \ref{fig:mbbBeamLsApproachesVoidEvolComparison} shows snapshots of the void domain to highlight the hole nucleation process of the LSO, SFC, and TFC beam designs in Fig. \ref{fig:mbbBeamLsAndDensOnlyAndCombAppMeshRef} using a mesh of size $h=1.0$. 
The evolution of the material removed by each approach is shown at $\mathcal{D}_{it}$=[5, 25, 50, 100].
The initial hole pattern in the LSO approach favors the removal of material along the external boundaries of the domain, which eventually results in a structure with internal shear-webs. In contrast, the holes nucleated in both the SFC and TFC are concentrated inside the design domain. Thus, in this particular case, the hole seeding promotes the formation of the shear-webs along the external boundaries.


\subsection{Example 3: bracket} \label{sec:Bracket3D}

Finally, the proposed hole nucleation approaches are applied to a complex engineering design problem.
The complexity of this design problem stems from the 3D geometry of the design domain and the multiple objective components and constraints.

\begin{figure}[h!]
	\center
	\includegraphics[width=0.4\linewidth]{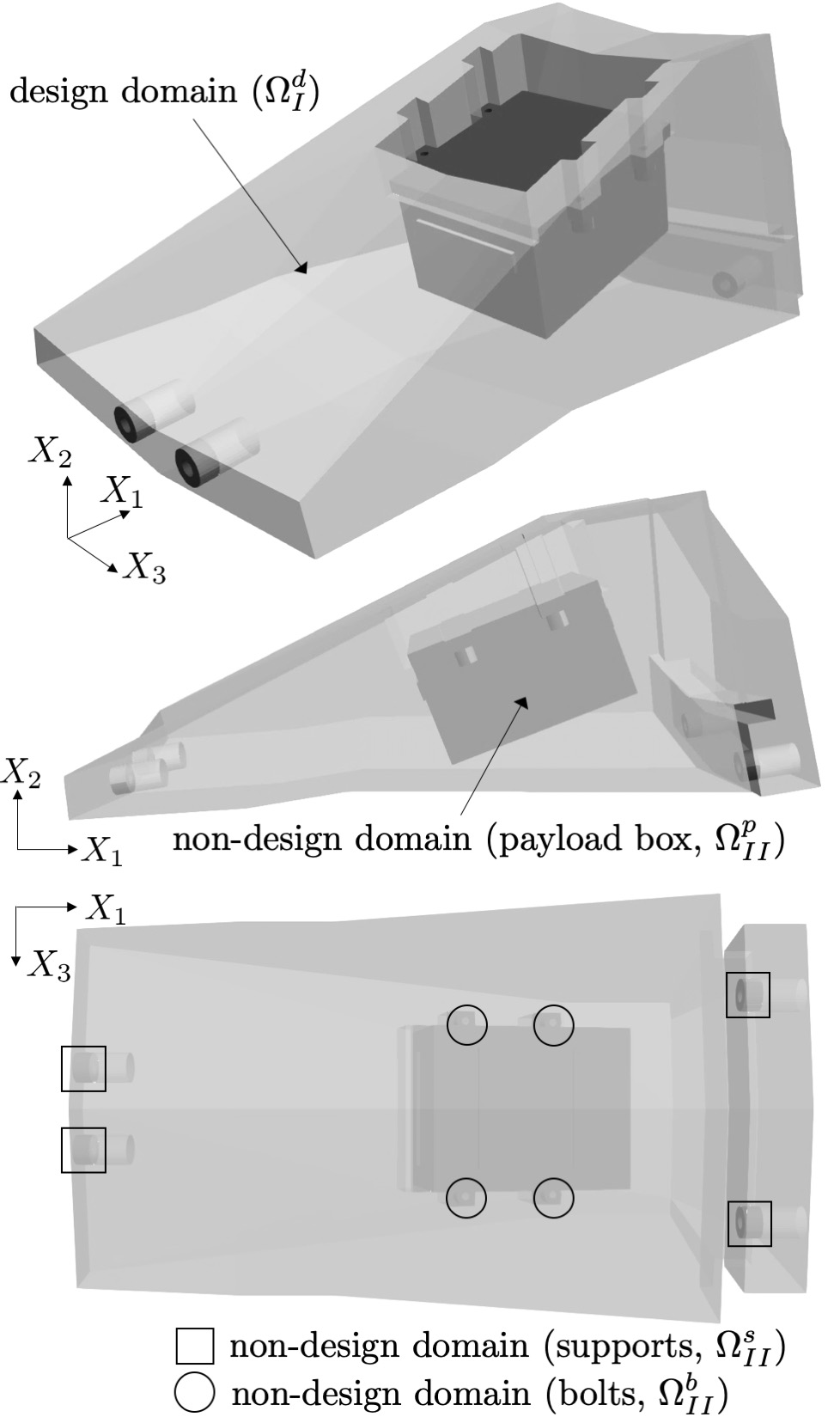}
	\caption{Bracket problem design and non-design domains in multiple views.}
	\label{fig:bracket3DProbSetup}
\end{figure}

The goal of this design problem is to find a structure that supports a payload given a set of supports and bolts for attaching the structure to the payload.
Fig. \ref{fig:bracket3DProbSetup} shows the design domain ($\Omega_{I}^{d}$) colored in light grey, while part of the non-design domains for the payload box ($\Omega_{II}^{p}$), the supports ($\Omega_{II}^{s}$), and bolts ($\Omega_{II}^{b}$) colored in dark grey.  Bolts and supports are modeled via hollow cylinders.
Side and bottom views are provided in addition to a 3D view to better visualize the constrained subdomains.
A uniform pressure load acts on the top surface of the payload box. The entire structure is subject to a body force in $X_2$ direction, representing an equivalent shock loading. 
The material properties and load conditions for this static problem can be found in Table \ref{tab:ex3BracketMatProps}.

\begin{table}[h]
	\caption{\label{tab:ex3BracketMatProps}Material properties and load conditions employed in bracket problem.}
	\center
	\renewcommand{\arraystretch}{1.2}
	\begin{tabular}{l|c}
		\hline
		Property                  & Value\\\hline
		Young's Modulus ($\Omega_{I}$, solid)  	 &  $E_S = 1.138$x$10^7$  [N/$cm^2$]\\
		Young's Modulus ($\Omega_{II}$, void)           &  $E_V = 1.138$x$10^{-1}$ [N/$cm^2$] \\
		Material Density ($\Omega_{I}$, solid)     	 		 &  $\theta_S = 4.43$x$10^{-5}$ [kg/$cm^3$]\\
		Material Density ($\Omega_{II}$, void)          		 &  $\theta_V = 0.0$ [kg/$cm^3$]\\
		Poisson Ratio ($\Omega_{I}$ and $\Omega_{II}$) &  $\nu_S = \nu_V$ = 0.342 [-] \\
		Pressure load ($T_{X_2}$) &  1.2x$10^4$ [N/$cm^2$] \\	
		Maximum stress ($\boldsymbol{\sigma}_{max}$) &  398.7 [N/$cm^2$] \\				
		\hline
	\end{tabular}
\end{table}

\begin{table}[h]
	\caption{\label{tab:Ex3DefContParams}Continuation parameters for the bracket problem.}
	\center
	\renewcommand{\arraystretch}{1.2}
	\begin{tabular}{l|c}
		\hline
		Parameter                  & Value\\\hline
		Continuation step size&  $\mathcal{D}_{st}=25$              \\
		Number of design iterations in continuation &  $\mathcal{D}_{c}=125$              \\		
		Maximum number of design iterations &  $\mathcal{D}_{max}=185$              \\			
		Continuation density threshold exponent &  $\eta_{\rho_{th}}=2.0$       \\
		Continuation density shift exponent &  $\eta_{\rho_{sh}}=2.0$              \\
		Initial density threshold &  $\rho^0_{th}=0.85 \rho_0$              \\
		Initial density shift &  $\rho^0_{sh}=0.1$              \\			
		\hline
	\end{tabular}
\end{table}

%

\begin{figure}[ht]\center
  \includegraphics[width=0.5\linewidth]{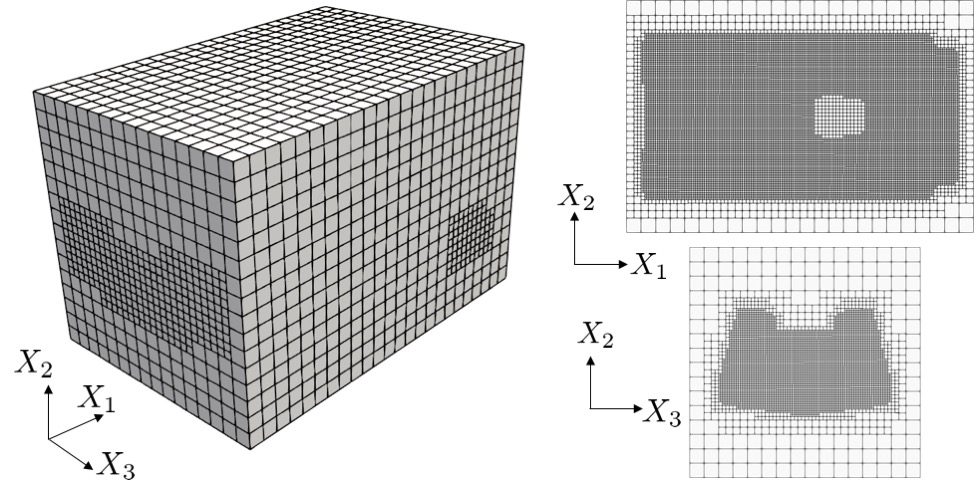}
\caption{3D and plane views of the locally refined mesh used in the bracket problem.}
\label{fig:bracketMeshViews}
\end{figure}

The objective of this optimization problem is to minimize the mass, $\mathcal{M}$, of the structure that supports the payload box.
The optimization problem formulation is the following:
\begin{equation}\label{eq:Ex3OptProbSetup}
\begin{aligned}
\underset{s}{\min}~z(\mathbf{s}) & = 
z_{sc} \bigg(
w_1~\mathcal{M} (\mathbf{s})  
+ 
w_2~\Psi (\mathbf{s},\mathbf{u(\mathbf{s})})/\Psi_0
+ 
w_3~P_{\hat\phi}(\mathbf{s}) \\
& + 
w_4~P_{Per}(\mathbf{s})
+
w_5~ P_{Reg}(\mathbf{s}) 
+ 
w_6~P_{\rho\phi}(\mathbf{s})
\bigg)
\\
s.t.: 
g_1 & = \frac{\mathcal{M} (\mathbf{s}) }{\Omega_{I}+\Omega_{II}} - \gamma_m\leq 0.
\\
g_2 & = w_{P_{\hat\phi}} ~ P_{\hat\phi}(\mathbf{s}) \leq 0.
\\
g_3 & = w_{P_{\tau}} ~  P_{\tau}({\mathbf{s},\mathbf{u(\mathbf{s})}}) \leq 0.
\end{aligned}
\end{equation}
The first term in the objective is the mass of the design domain to be minimized. 
The second term is a strain energy component added to prevent an overly aggressive removal of mass early in the optimization process and promote a sufficiently stiff structure.
The penalty term, $P_{\hat\phi}$, is formulated as follows:
\begin{equation}\label{eq:LsBracketNonDomPen}
\begin{aligned}
P_{\hat\phi} & = 
\frac{\displaystyle\int_{\Omega_D} ( \phi(\mathbf{X}) - \hat\phi(\mathbf{X}) )^2 dV }{\displaystyle\int_{\Gamma_D} dA},
\end{aligned}
\end{equation}
with the field $\hat\phi(\mathbf{X})$ defined as:
\begin{equation}\label{eq:tildeLsFuncBracket}
\begin{aligned}
\hat\phi(\mathbf{X}) & = 
\begin{cases}
\phi_{low}, & \forall~\boldsymbol X \in \Omega_D / (\Omega^{d}_{I} \cup \Omega^{p}_{II} \cup \Omega^{s}_{II} \cup \Omega^{b}_{II}), \\
\phi_{low}, & \forall~\boldsymbol X \in \Gamma_{(\Omega^{d}_{I} \cap \Omega^{p}_{II}) / \Omega^{b}_{II}  }, \\
\phi_{up}, & \forall~\boldsymbol X \in (\Omega^{p}_{II} \cup \Omega^{s}_{II} \cup \Omega^{b}_{II}), \\
\phi(\mathbf{X}), & otherwise;
\end{cases}
\end{aligned}
\end{equation}
based on the previously defined subdomains highlighted in Fig. \ref{fig:bracket3DProbSetup}.
Through the field defined in Eq. \ref{eq:tildeLsFuncBracket}, the LS penalization in Eq. \ref{eq:LsBracketNonDomPen} (i) forces the support structure to stay within the design domain, (ii) prevents contact between the payload box and the support structure except for the bolts, and (iii) avoids altering the LSF in the payload box, attachment bolts and the supports.
The remaining penalty terms in the objective ($P_{Per}$, $P_{Reg}$, and $P_{\rho\phi}$) are defined in Section \ref{subsec:TFCOptProbForm}.
Similar to the beam example, the perimeter weight, $w_4$, is updated at every continuation step from 0.0001 to 0.01, using Eq. \ref{eq:expFuncContInc} with an exponent of 3.0.

In addition to considering the mass in the objective, a mass constraint ($g_1$, with $\gamma_m=0.30$) is imposed to represent an upper limit on the mass. 
The stress constraint, $g_3$, enforces the volume of solid in which the stress exceeds the maximum von Mises stress, $\sigma_{VM}^{max}$, to zero. The stress penalty,
\begin{equation}\label{eq:stressConstForm}
\begin{aligned}
P_{\tau}
=
\int_{\Gamma_{D}} \hat\tau(\mathbf{X}) dV,
\end{aligned}
\end{equation}
is an integral over the volume of the scalar stress field, $\hat\tau(\mathbf{X})$, defined as:
\begin{equation}\label{eq:stressFuncForm}
\begin{aligned}
\hat\tau(\mathbf{X})
=
\begin{cases}
\left[  (\tau- \sigma_{VM}^{max})^2 + \xi^2_{\tau} \right]^{1/2} - \xi_{\tau}, & \forall~\tau- \sigma_{VM}^{max}>0,\\
0, & \forall~\tau- \sigma_{VM}^{max} \leq 0.
\end{cases} 
\end{aligned}
\end{equation}
The parameter $ \xi_{\tau}$ is set to 0.1, and $\sigma_{VM}^{max}$ is given by the yield stress of Ti-6Al-4V, reduced by a safety factor of 2.0 (see Table \ref{tab:ex3BracketMatProps}).
This last constraint remains inactive initially but controls the mass removal effect towards the final stages of the optimization process.

The weights employed in this problem in the objective components and constraints are $w_1=5.0$, $w2 = 0.005$, $w_3=w_{P_{\hat\phi}}=5000.0$, $w_4=0.0001$, $w_5=0.01$, $w_6=10.0$, and $w_{P_{\tau}}=1.0\times10^4$.
The objective scaling parameter is set to $z_{sc} = 10.0$. 

\begin{figure*}[h!]
	\center
	\includegraphics[width=1.0\linewidth]{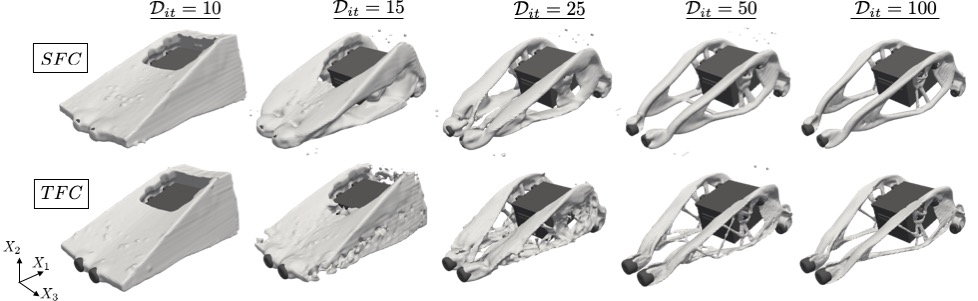}
	\caption{Snapshots of the bracket design problem at iterations $\mathcal{D}_{it}$= [10, 15, 25, 50, 100].}
	\label{fig:bracketEvolStrEnAngMassWithDesigns}
\end{figure*}

The design optimization problem is embedded in a computational domain of 15.24$\times$10.16$\times$10.16 $cm^3$, and solved using the box mesh in Fig. \ref{fig:bracketMeshViews}. 
The mesh is locally refined within the design domain to reduce the computational cost.
A 3D view of the entire mesh and two plane views at the center of the mesh in the $X_1-X_2$ and $X_2-X_3$ planes are shown on the top right and bottom right sides, respectively.
The design domain, supports and bolts are immersed into this mesh. 
Note that the mesh transitions occur outside the design domain $\Omega_I^d$; hence, no special treatment of hanging nodes is needed in the optimization framework employed. 

The initial density and the SIMP exponent are set to $\rho_0=0.2$ and $\beta{\rho}=2.0$, respectively. The continuation parameters for the density shift and density threshold can be found in Table \ref{tab:Ex3DefContParams}.

\begin{figure}[h!]
	\center
	\includegraphics[width=0.5\linewidth]{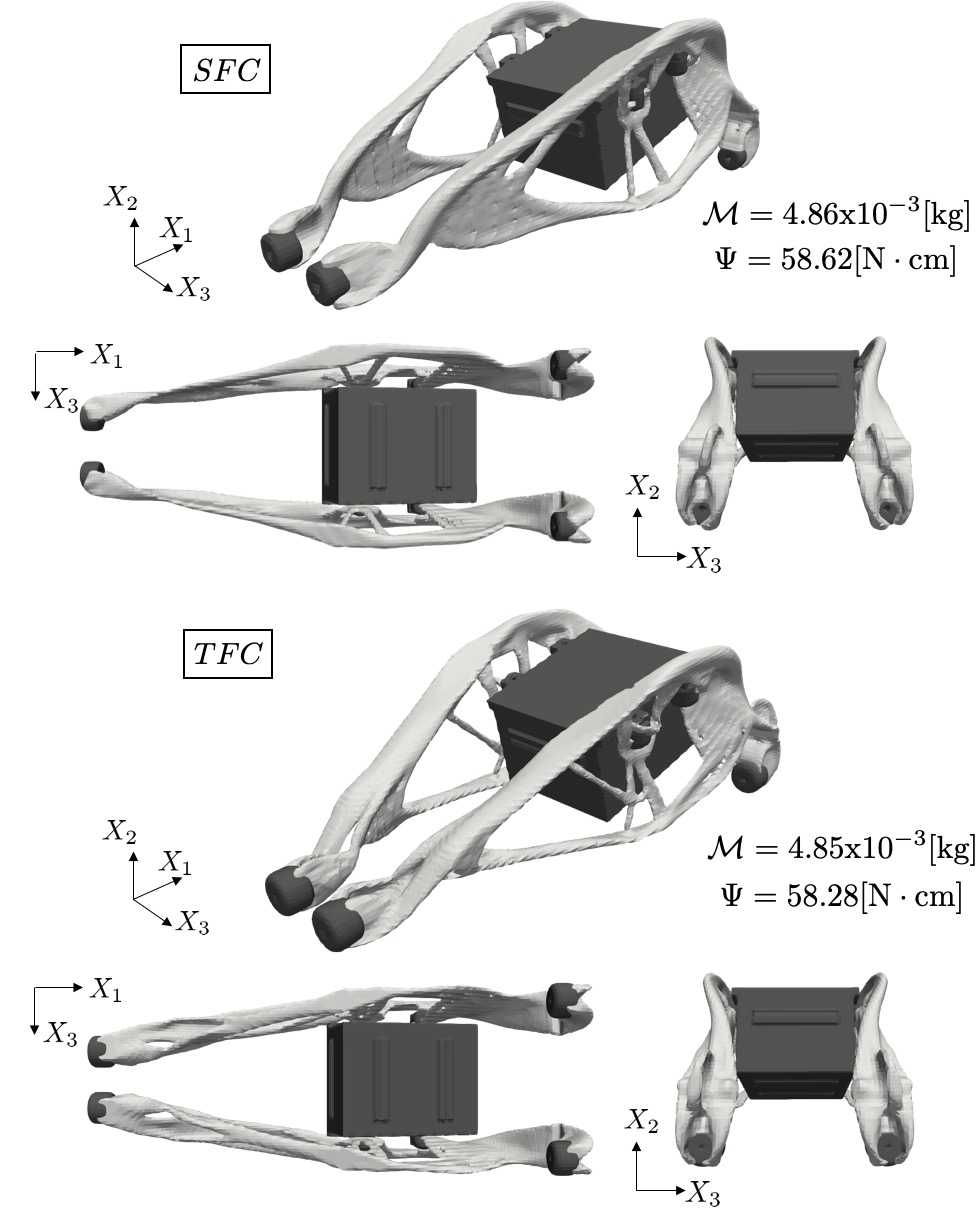}
	\caption{Final bracket designs for both the SFC (top) and TFC (bottom) approaches.}
	\label{fig:bracket3DFinalDesigns}
\end{figure}

Snapshots of the designs for both coupling strategies are shown in Fig. \ref{fig:bracketEvolStrEnAngMassWithDesigns}. The structure colored in light grey evolves during the optimization process while the non-design domains in dark grey remain fixed.
The TFC approach favors removing material internally, while the SFC approach nucleates less holes and evolves the design along the solid-void interface instead. This is observed predominantly for $\mathcal{D}_{it} = [15,25]$ in  Fig. \ref{fig:bracketEvolStrEnAngMassWithDesigns}.
In both cases, most of the topological changes occur within the first 30 design iterations consequence of nearly all holes being nucleated early in the optimization process. 
Well-defined geometries that satisfy the geometric and stress constraints are attained in the final designs, as shown in Fig. \ref{fig:bracket3DFinalDesigns}. 
The SFC scheme converges to a design with more thin-walled features in the front section of the domain, while the TFC final design shows more truss-like features.
Despite the topological differences, both approaches converge to final designs with similar performances.

\begin{figure*}[h]
	\center
	\includegraphics[width=1.0\linewidth]{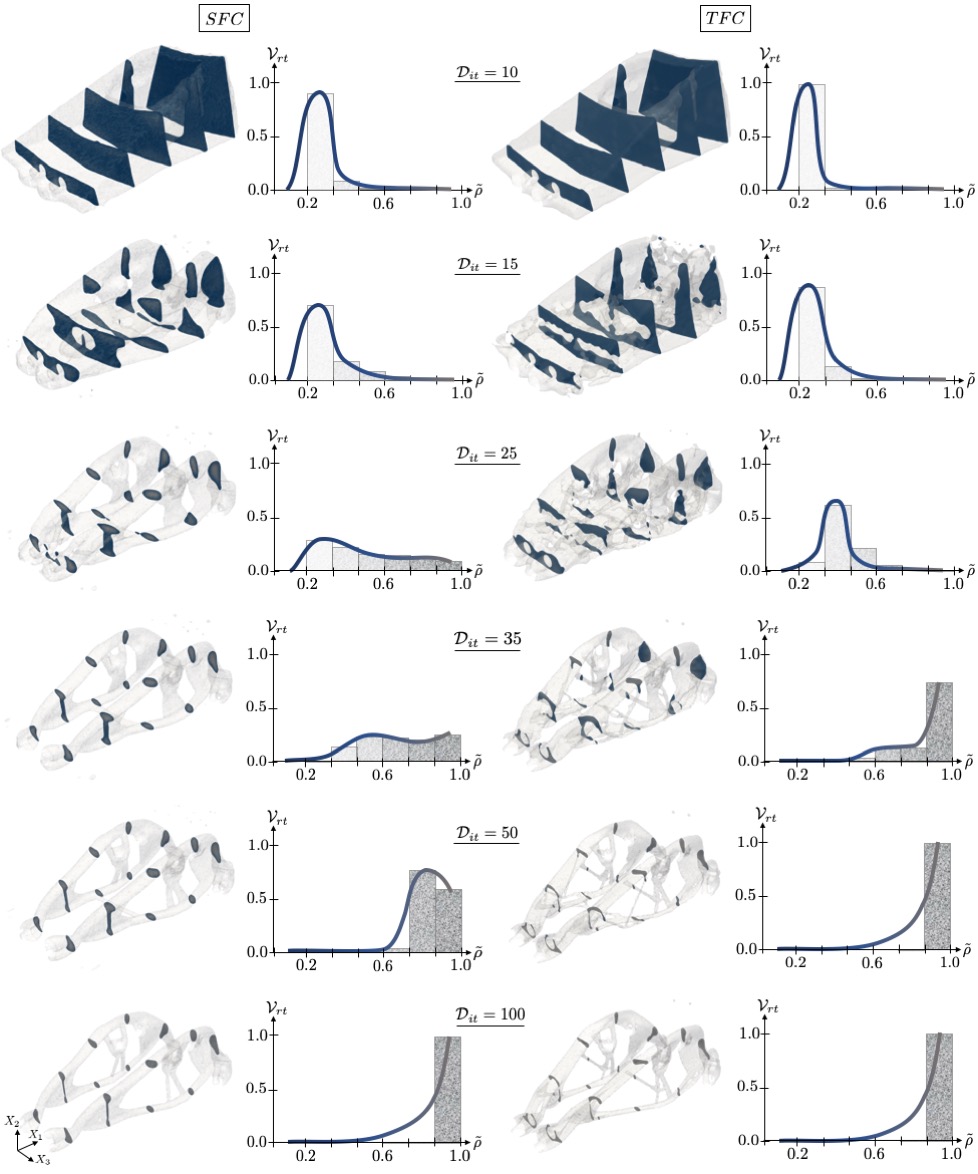}
	\caption{Snapshots of bracket problem at $\mathcal{D}_{it}$=[10, 15, 25, 35, 50, 100] together with histograms of their density distributions using SFC (left) and TFC (right) schemes are shown.}
	\label{fig:bracket3DVertCutsWithDensFieldAndHist}
\end{figure*}

Cuts of the design domains of the bracket at $\mathcal{D}_{it}$=[10, 15, 25, 35, 50, 100] are shown for both approaches in Fig. \ref{fig:bracket3DVertCutsWithDensFieldAndHist}. The $X_3-X_2$ plane views at [4.24, 6.73, 9.22, 11.73, 14.73] $cm$ away from the front boundary are colored by the shifted density field.
Each design is accompanied by a histogram of the density distribution.
The histograms were constructed by computing a volume ratio, $V_{rt}$, in seven density intervals. 
The parameter $V_{rt}$ is defined as the ratio between the volume of material contained on each density interval and the total material volume of the design domain.
Despite starting from the same unbiased initial shifted density distribution, clear differences in their evolution are observed between the single and two-field approaches. 
The histograms illustrate more disperse shifted density distributions for the SFC throughout the hole nucleation process. As a consequence, the SFC approach shows less pronounced topological changes at early stages of the optimization process.
Although in both cases a [0-1] density distribution is achieved due to the density shift, defining the density field in TFC approach by an independent set of variables allows for faster convergence. 
Intermediate densities present in the vicinity of the material interface by construction in the SFC approach partially restrict, and thus decelerate, the evolution of the design.

Overall, the results show that the proposed SFC and TFC LS TO approaches are well suited to tackle complex structural optimization problems with stress constraints. 

\section{Conclusions and Future Work} \label{sec:Concl}

Two LS TO approaches that use a density field to nucleate holes during the optimization process and accelerate convergence were developed and studied.
LS and density fields were coupled via (i) a single vector of abstract optimization variables (SFC); and (ii) two vectors of optimization variables that govern them independently through a penalty term in the objective (TFC). 
In both cases no initial seeding is required, and the LS and density fields are optimized simultaneously.
Contrary to classical density-based methods, which require sufficient penalization of intermediate densities to converge to a [0-1] material distribution, here intermediate densities are eliminated by introducing a density shift. A full material density distribution is achieved by gradually shifting the densities through a continuation scheme.

Benchmark numerical examples in 2D and 3D, as well as a geometrically complex engineering application demonstrated that the proposed approaches are overall robust with respect to various algorithmic parameters. 
We also observed some influence of the density penalization on the optimized design.
These coupling strategies can be extended to incorporate both multiple LSs and density interpolation schemes using either explicit or implicit LS TO approaches.

In future work, the proposed approach should be studied for problems with different objectives and constraints, and for problems involving different physics models.
Future work could also investigate enhanced formulation of the SFC and TFC approaches to gain feature size control and robustness against shape imperfections.

\section*{Acknowledgements}

All authors acknowledge the support of the National Science Foundation (CMMI-1463287). The third author acknowledges the support of the Defense Advanced Research Projects Agency (DARPA) under the TRADES program (agreement HR0011-17-2-0022). The opinions and conclusions presented in this paper are those of the authors and do not necessarily reflect the views of the sponsoring organizations.



\end{document}